# ESTIMATION AND CONFIDENCE SETS FOR SPARSE NORMAL MIXTURES

By T. Tony Cai,[1] Jiashun Jin[2] and Mark G. Low[1]

*University of Pennsylvania, Purdue University and University of Pennsylvania*

For high dimensional statistical models, researchers have begun to focus on situations which can be described as having relatively few *moderately large* coefficients. Such situations lead to some very subtle statistical problems. In particular, Ingster and Donoho and Jin have considered a sparse normal means testing problem, in which they described the precise demarcation or *detection boundary*. Meinshausen and Rice have shown that it is even possible to estimate consistently the fraction of nonzero coordinates on a subset of the detectable region, but leave unanswered the question of exactly in which parts of the detectable region consistent estimation is possible.

In the present paper we develop a new approach for estimating the fraction of nonzero means for problems where the nonzero means are moderately large. We show that the detection region described by Ingster and Donoho and Jin turns out to be the region where it is possible to consistently estimate the expected fraction of nonzero coordinates. This theory is developed further and minimax rates of convergence are derived. A procedure is constructed which attains the optimal rate of convergence in this setting. Furthermore, the procedure also provides an honest lower bound for confidence intervals while minimizing the expected length of such an interval. Simulations are used to enable comparison with the work of Meinshausen and Rice, where a procedure is given but where rates of convergence have not been discussed. Extensions to more general Gaussian mixture models are also given.

**1. Introduction.** In many statistical applications such as analysis of microarray data, signal recovery and functional magnetic resonance imaging

---

Received November 2005; revised December 2006.

[1]Supported in part by NSF Grants DMS-03-06576 and DMS-06-04954.

[2]Supported in part by NSF Grant DMS-05-05423.

(fMRI), the focus is often on identifying and estimating a relatively few significant components from a high dimensional vector. In such applications, models which allow a parsimonious representation have important advantages, since effective procedures can often be developed based on relatively simple testing and estimation principles. For example, in signal and image recovery, wavelet thresholding is an effective approach for recovering noisy signals since wavelet expansions of common functions often lead to a sparse representation; the quality of the recovery depends only on the large coefficients; the "small" coefficients have relatively little effect on the quality of the reconstruction, and thresholding rules are effective in identifying and estimating the large coefficients. Likewise, in problems of multiple comparison where only a very small fraction of hypotheses are false, the false discovery rate (FDR) approach introduced by Benjamini and Hochberg [1] is an effective tool for identifying those false hypotheses.

In these problems, the focus is on discovering large components. However, recently there has been a shift of attention toward problems which involve identifying or estimating "moderately" large components. Such terms cannot be isolated or detected with high probability individually. However it is possible to detect the presence of a collection of such "moderate" terms. For multiple comparison problems where there are a large number of tests to be performed, it may not be possible to identify the particular false hypotheses, although it is possible to discover the fraction of the false null hypotheses. For example, Meinshausen and Rice [14] discuss the Taiwanese–American Occultation Survey, where it is difficult to tell whether an occultation has occurred for a particular star at a particular time, but it is possible to estimate the fraction of occultations that have occurred over a period of time. In this setting, it is not possible to perform individual tests with high precision, but it is possible to estimate the fraction of false nulls. Other examples include the analysis of Comparative Genomic Hybridization (CGH) lung cancer data [11], microarray breast cancer data [6, 10] and Single Nucleotide Polymorphism (SNP) data on Parkinson disease [13].

For such applications where there are relatively few nonzero components, it is natural to develop the theory with a random effects model; see, for example, Efron [6], Meinshausen and Rice [14] and Genovese and Wasserman [7]. Consider $n$ independent observations from a Gaussian mixture model,

$$(1.1) \qquad X_i = \mu_i + z_i, \qquad z_i \stackrel{\text{i.i.d.}}{\sim} N(0, 1), \qquad 1 \le i \le n,$$

where $\mu_i$ are the random effects with $P(\mu_i = 0) = 1 - \varepsilon_n$, and given $\mu_i \ne 0$, $\mu_i \sim H$ for some distribution $H$. Equivalently we may write

$$(1.2) \qquad X_i \stackrel{\text{i.i.d.}}{\sim} (1 - \varepsilon_n) N(0, 1) + \varepsilon_n G, \qquad 1 \le i \le n,$$

where $G$ is the convolution between $H$ and a standard Gaussian distribution. In these models, the problem of estimating the fraction of nonzero terms corresponds to estimating the parameter $\varepsilon_n$, and we are particularly interested in the case where the signal is sparse and the nonzero terms $\mu_i$ are "moderately" large (i.e., $\varepsilon_n$ is small and $|\mu_i| < \sqrt{2 \log n}$). This general problem appears to be of fundamental importance.

The development of useful estimates of $\varepsilon_n$ along with the corresponding statistical analysis appears to pose many challenges. In fact this theory is already quite involved even in the apparently simple special case where $H$ is concentrated at a single point $\mu_n$; here $\mu_n$ depends on $n$ but not on $i$. In this case (1.2) becomes a two-point mixture model,

$$(1.3) \qquad X_i \overset{\text{i.i.d.}}{\sim} (1 - \varepsilon_n) N(0, 1) + \varepsilon_n N(\mu_n, 1), \qquad 1 \le i \le n.$$

In such a setting, the problem of testing the null hypothesis $H_0 : \varepsilon_n = 0$ against the alternative $H_a : \varepsilon_n > 0$ was first studied in detail in Ingster [8], where $(\varepsilon_n, \mu_n)$ are assumed to be known (see also [9]). Ingster showed that this apparently simple testing problem contains a surprisingly rich theory even though the optimal test is clearly the likelihood ratio test. Donoho and Jin [5] extended this work to the case of unknown $(\varepsilon_n, \mu_n)$. It was shown that the interesting range for $(\varepsilon_n, \mu_n)$ corresponds to a relatively "small" $\varepsilon_n$ and a "moderately" large $\mu_n$. A detection boundary was developed which separates the possible pairs $(\varepsilon_n, \mu_n)$ into two regions, the *detectable* region and the *undetectable* region. When $(\varepsilon_n, \mu_n)$ belongs to the interior of the undetectable region, the null and alternative hypotheses merge asymptotically and no test could successfully separate them. When $(\varepsilon_n, \mu_n)$ belongs to the interior of the detectable region, the null and alternative hypotheses separate asymptotically.

Although the theory of testing the above null hypothesis is closely related to the estimation problem we are considering, it does not automatically yield estimates of $\varepsilon_n$. In fact, the problem of estimating $\varepsilon_n$ appears to contain further challenges which are not present in the above testing problem. Even the theory for consistent estimation of $\varepsilon_n$ recently studied in Meinshausen and Rice [14] is quite complicated. Meinshausen and Rice [14] gave an estimate of $\varepsilon_n$ and showed it to be consistent on a subset of the detectable region. They pointed out that "it is clear that it is somewhat easier to test for the global null hypothesis than to estimate the proportion," leaving the following question unanswered: what is the precise region over which consistent estimation of $\varepsilon_n$ is possible?

There are two primary goals of the present paper. The first is to develop in detail the theory for estimating $\varepsilon_n$ in the two-point Gaussian mixture model. The theory given in the present paper goes beyond consistent estimation, and focuses on the development of procedures which have good mean squared error performance. Minimax rates of convergence are shown to depend on the magnitude of both $\mu_n$ and $\varepsilon_n$; upper and lower bounds for the minimax mean squared error are given, which differ only by logarithmic factors; estimates of $\varepsilon_n$ which adapt to the unknown $\mu_n$ and $\varepsilon_n$ are also given. These results make precise how accurately $\varepsilon_n$ can be estimated in such a model. In particular, we show that it is possible to estimate $\varepsilon_n$ consistently whenever $(\varepsilon_n, \mu_n)$ is in the detectable region; and although the estimation problem is in some sense technically more challenging than the testing problem, the estimable region and detectable region actually coincide.

The other major goal of the present paper is to show that the theory developed for the two-point mixture model leads to a one-sided confidence interval for $\varepsilon_n$, which has guaranteed coverage probability not only for the two-point mixture model, but also over the mixture model (1.1) assuming only that $H > 0$. In this general one-sided Gaussian mixture model, as noted in a similar context by Meinshausen and Rice [14], the upper bound for $\varepsilon_n$ must always be equal to 1: the possibility that $\varepsilon_n = 1$ can never be ruled out because the nonzero $\mu_i$ can be arbitrarily close to zero. For example, asymptotically it is impossible to tell whether all the $\mu_i$ are zero or all of them are equal to, say, $10^{-n}$. On the other hand, if many "large" values of $X_i$ are observed it is possible to give useful lower bounds on the value of $\varepsilon_n$. This is therefore an example of a situation where only one-sided inference is possible; a nontrivial lower bound for $\varepsilon_n$ can be given but not a useful upper bound. See Donoho [4] for other examples and a general discussion of problems of one-sided inference. In such a setting, a natural goal is to provide a one-sided confidence interval for the parameter of interest, which both has a guaranteed coverage probability and is also "close" to the unknown parameter. We show that such a one-sided confidence interval can be built by using the theory developed for the two-point model.

The paper is organized as follows. We start in Section 2 with the two-point mixture model. As mentioned earlier, this model has been the focus of recent attention both for testing the null hypothesis that $\varepsilon_n = 0$ and for consistent estimation of $\varepsilon_n$. These results are briefly reviewed and then a new family of estimators for $\varepsilon_n$ is introduced. A detailed analysis of these estimators requires precise bounds on the probability of over-estimating $\varepsilon_n$, which can be given in terms of the probability that a particular confidence band covers the true distribution function. Section 3 is devoted to giving accurate upper bounds of this probability. In Section 4 we consider the implication of these results for estimating $\varepsilon_n$ under mean squared error. Section 5 is devoted to the theory of one-sided confidence intervals over all one-sided mixture models. Section 6 connects the results of the previous sections to that of consistent estimation of $\varepsilon_n$, where comparisons to the work of Meinshausen and Rice [14] are also made. While the above theory is asymptotic, the discussion is continued in Section 7, where simulations show that the procedure performs well in settings similar to those considered by Meinshausen and Rice. Proofs are given in Section 8.

**2. Estimation of $\varepsilon_n$ in the two-point mixture model.** In this section we focus on estimating the fraction $\varepsilon_n$ under the two-point mixture model,

$$(2.1) \qquad X_i \overset{\text{i.i.d.}}{\sim} (1 - \varepsilon_n)N(0,1) + \varepsilon_n N(\mu_n, 1), \qquad 1 \le i \le n.$$

As mentioned in the Introduction, the problems of testing the null hypothesis that $\varepsilon_n = 0$ and estimating $\varepsilon_n$ consistently in the sense that $P\{|\frac{\hat{\varepsilon}_n}{\varepsilon_n} - 1| > \delta\} \to 0$ for all $\delta > 0$ have been considered. These results are briefly reviewed in Section 2.1 so as to help clarify the goal of the present work. A new family of estimators is then introduced in Section 2.2. Later sections show how to select from this family of estimators those which have good mean squared error performance, and those which provide a lower end point for a one-sided confidence interval with a given guaranteed coverage probability.

2.1. *Review of testing and consistency results.* Ingster [8] and Donoho and Jin [5] studied the problem of testing the null hypothesis that $\varepsilon_n = 0$. It was shown that the interesting cases correspond to choices of $\varepsilon_n$ and $\mu_n$ where $(\varepsilon_n, \mu_n)$ are calibrated with a pair of parameters $(r, \beta)$: $\varepsilon_n = n^{-\beta}$ and $\mu_n = \sqrt{2r \log n}$, where $1/2 < \beta < 1$ and $0 < r < 1$. Under this calibration it was shown that there is a detection boundary which separates the testing problem into two regions. Set

$$(2.2) \qquad \rho^*(\beta) = \begin{cases} \beta - \frac{1}{2}, & 1/2 < \beta \le 3/4, \\ (1 - \sqrt{1 - \beta})^2, & 3/4 < \beta < 1. \end{cases}$$

In the $\beta$-$r$ plane, we call the curve $r = \rho^*(\beta)$ the *detection boundary* [5, 8, 9] associated with this hypothesis testing problem. The detection boundary separates the $\beta$-$r$ plane into two regions: the *detectable region* and the *undetectable region.* When $(\beta, r)$ belongs to the interior of the undetectable region, the sum of Type I and Type II errors for testing the null hypothesis that $\varepsilon_n = 0$ against the alternative $(\varepsilon_n = n^{-\beta}, \mu_n = \sqrt{2r \log n})$ must tend to 1. Hence no test can asymptotically distinguish the two hypotheses. On the other hand when $(\beta, r)$ belongs to the interior of the detectable region, there are tests for which both Type I and Type II errors tend to zero and thus the hypotheses can be separated asymptotically. These two regions are illustrated in Figure 2, where a third region—the classifiable region—is also displayed. When $(\beta, r)$ belongs to the interior of the classifiable region, it is not only possible to reliably tell that $\varepsilon_n > 0$, but also to separate the observations into signal and noise.

It should be stressed that this testing theory does not yield an effective strategy for estimating $\varepsilon_n$, though it does provide a benchmark for a theory of consistent estimation. Important progress in this direction has recently been made in Meinshausen and Rice [14], where an estimator of $\varepsilon_n$ was constructed and shown to be consistent if $r > 2\beta - 1$. This estimator is however inconsistent when $r < 2\beta - 1$. Note here that the separating line $r = 2\beta - 1$ always falls above the detection boundary. See Figure 2. The work of Meinshausen and Rice leaves unclear the question of whether consistent estimation of $\varepsilon_n$ is possible over the entire detectable region. Of course, in the undetectable region no estimator can be consistent, as any consistent estimator immediately gives a reliable way for testing $\varepsilon_n = 0$.

2.2. *A family of estimators.* The previous section outlined the theory developed to date for estimating $\varepsilon_n$ in the two-point Gaussian mixture model (1.3). The goal of the present paper is to develop a much more precise estimation theory both for one-sided confidence intervals as well as for mean squared error. A large part of this theory relies on the construction of a family of easily implementable procedures along with an analysis of particular estimators chosen from this family of estimators. The present section focuses on providing a detailed description of the construction of this family of estimators. Later in Sections 4 and 5 we will show how to choose particular members of this family to yield near optimal mean squared error estimates and one-sided confidence intervals.

The basic idea underlying the general construction given here relies on the following representation for $\varepsilon_n$. Throughout the paper we shall denote by $\phi$ and $\Phi$, respectively, the density and cumulative distribution function (c.d.f.) of a standard normal distribution. Suppose that instead of observing the data (2.1), one can observe directly the underlying c.d.f. $F(t) \equiv (1 - \varepsilon_n)\Phi(t) + \varepsilon_n\Phi(t - \mu_n)$ at just two points, say $\tau$ and $\tau'$ with $0 \le \tau < \tau'$. Then the values of $\varepsilon_n$ and $\mu_n$ can be determined precisely as follows. Set

$$(2.3) \qquad D(\mu; \tau, \tau') = [\Phi(\tau) - \Phi(\tau - \mu)] / [\Phi(\tau') - \Phi(\tau' - \mu)].$$

Lemma 8.1 in [2] shows that $D(\cdot; \tau, \tau')$ is strictly decreasing in $\mu > 0$ for any $\tau < \tau'$. The parameters $\varepsilon_n$ and $\mu_n$ are then uniquely determined by

$$(2.4) \qquad \varepsilon_n = \frac{\Phi(\tau) - F(\tau)}{\Phi(\tau) - \Phi(\tau - \mu_n)} \quad \text{and} \quad D(\mu_n; \tau, \tau') = \frac{\Phi(\tau) - F(\tau)}{\Phi(\tau') - F(\tau')}.$$

It is easy to check that for $\tau < \tau'$,

$$\inf_{\mu > 0} D(\mu; \tau, \tau') \equiv \frac{\Phi(\tau)}{\Phi(\tau')} < \frac{\Phi(\tau) - F(\tau)}{\Phi(\tau') - F(\tau')} < \sup_{\mu > 0} D(\mu; \tau, \tau') \equiv \frac{\phi(\tau)}{\phi(\tau')},$$

so by the monotonicity of $D(\cdot; \tau, \tau')$, we can first solve for $\mu_n$ from the right-hand side equation in (2.4), and then plug this $\mu_n$ into the left-hand side equation in (2.4) for $\varepsilon_n$.

In principle estimates of $\mu_n$ and $\varepsilon_n$ can be given by replacing $F(\tau)$ and $F(\tau')$ by their usual empirical estimates. Unfortunately, this simple approach does not work well since the performance of the resulting estimate depends critically on the choice of $\tau$ and $\tau'$. For most choices of $\tau$ and $\tau'$ the resulting estimate is not a good estimate of $\varepsilon_n$ in terms of mean squared error, although it is often consistent. Moreover, although there are particular pairs for which the resulting estimator does perform well, it is difficult to select the optimal pair of $\tau$ and $\tau'$ since the optimal choice depends critically on the unknown parameters $\varepsilon_n$ and $\mu_n$. It is however worth noting that for the situations considered here the optimal choices of $\tau$ and $\tau'$ always satisfy $0 \leq \tau < \tau' \leq \sqrt{2 \log n}$.

The key to the construction given below is that, instead of using the usual empirical c.d.f. as estimates of $F(\tau)$ and $F(\tau')$, we use slightly biased estimates of these quantities to yield an estimate of $\varepsilon_n$ which is with high probability smaller than the true $\varepsilon_n$. It is in fact important to do this over a large collection of $\tau$ and $\tau'$ so that the entire collection of estimates is simultaneously smaller than $\varepsilon_n$ with large probability. It then follows that the maximum of these estimates is also smaller than $\varepsilon_n$ with this same high probability. This resulting estimate is just one member of our final family of estimates; other members of this family are found by adjusting the probability that the initial collection of estimators underestimates $\varepsilon_n$. The details of this construction are given below.

First note that underestimates of $\varepsilon_n$ can be obtained by overestimating $F(\tau)$ and underestimating $F(\tau')$. More specifically, suppose that $F^+(\tau) \geq F(\tau)$ and $F^-(\tau') \leq F(\tau')$. Then there are two cases depending on whether or not the following holds:

$$\frac{\Phi(\tau)}{\Phi(\tau')} \leq \frac{\Phi(\tau) - F^+(\tau)}{\Phi(\tau') - F^-(\tau')} \leq \frac{\phi(\tau)}{\phi(\tau')}.$$

If it does not hold, then the equation does not give a good estimate for $\mu_n$ and we take 0 to be an estimate for $\varepsilon_n$. If it does hold, then we can use (2.4) to estimate $\mu_n$ by simply replacing $F(\tau)$ and $F(\tau')$ by $F^+(\tau)$ and $F^-(\tau')$, respectively. Call this estimate $\hat{\mu}_n$ and note that $\hat{\mu}_n \geq \mu_n$. It then immediately follows that the solution to the first equation in (2.4) with $\hat{\mu}_n$ replacing $\mu_n$ yields an estimate $\hat{\varepsilon}_n$ of $\varepsilon_n$ for which $\hat{\varepsilon}_n \leq \varepsilon_n$. A final estimator is then created by taking the maximum of these estimators.

Of course in practice we do not create estimators which always overestimate $F(\tau)$ and underestimate $F(\tau')$, as there is also another goal, namely that these estimates are also close to $F(\tau)$ and $F(\tau')$. To reconcile these goals it is convenient to first construct a confidence envelope for $F(t)$. First fix a value $a_n$ and solve for $F(t)$: $\sqrt{n} \frac{|F_n(t) - F(t)|}{\sqrt{F(t)(1 - F(t))}} = a_n$, where $F_n$ is the usual empirical c.d.f. The result is a pair of functions $F_{a_n}^{\pm}(t)$,

$$(2.5) \quad F_{a_n}^{\pm}(t) = \frac{2F_n(t) + a_n^2/n \pm \sqrt{a_n^2/n + (4F_n(t) - 4F_n^2(t)) \cdot (a_n/\sqrt{n})}}{2(1 + a_n^2/n)}.$$

Note $F_n^-(t) \leq F(t) \leq F_n^+(t)$ if and only if $\sqrt{n}\frac{|F_n(t)-F(t)|}{\sqrt{F(t)(1-F(t))}} \leq a_n$. So for any $S_n \subseteq (-\infty, \infty)$ if we take $a_n$ to be the $\alpha$-upper percentile of $\sup_{t \in S_n}\{\sqrt{n}\frac{|F_n(t)-F(t)|}{\sqrt{F(t)(1-F(t))}}\}$, then $F_n^\pm(t)$ together give a simultaneous confidence envelope for $F(t)$ for all $t \in S_n$. For each $a_n$ the confidence envelope can then be used to construct a collection of estimators as follows. Pick equally spaced grid points over the interval $[0, \sqrt{2\log n}]$: $t_j = (j-1)/\sqrt{2\log n}$, $1 \leq j \leq 2\log(n) + 1$. For a pair of adjacent points $t_j$ and $t_{j+1}$ in the grid let $\hat{\mu}_{a_n}^{(j)} = \hat{\mu}_{a_n}^{(j)}(t_j, t_{j+1}; n, \Phi, F^+, F^-)$ be the solution of the equation

$$(2.6) \qquad D(\mu; t_j, t_{j+1}) = \frac{\Phi(t_j) - F_{a_n}^+(t_j)}{\Phi(t_{j+1}) - F_{a_n}^-(t_{j+1})},$$

when such a solution exists. If there is no solution set $\hat{\varepsilon}_j = 0$. Note that if a solution exists and $F$ lies in the confidence envelope (2.5), then $F_{a_n}^+(t_j) \geq F(t_j)$ and $F_{a_n}^-(t_{j+1}) \leq F(t_{j+1})$ and hence $\hat{\mu}_{a_n}^{(j)} \geq \mu_n$. It then also follows that

$$(2.7) \qquad \hat{\varepsilon}_{a_n}^{(j)} = \frac{\Phi(t_j) - F_{a_n}^+(t_j)}{\Phi(t_j) - \Phi(t_j - \hat{\mu}_j)}$$

satisfies $\hat{\varepsilon}_{a_n}^{(j)} \leq \varepsilon$. The final estimator $\hat{\varepsilon}_{a_n}^*$ is defined by taking the maximum of $\{\varepsilon_{a_n}^{(j)}\}$:

$$(2.8) \qquad \hat{\varepsilon}_{a_n}^* \equiv \max_{1 \leq j \leq 2\log n} \hat{\varepsilon}_{a_n}^{(j)}.$$

**3. Evaluating the probability of underestimation.** A family of estimators depending on $a_n$ was introduced in Section 2 in terms of a confidence envelope. A detailed analysis of these estimators depends critically on upper bounding the probability of overestimating $\varepsilon_n$. Note that $\hat{\varepsilon}_{a_n}^*$ underestimates $\varepsilon_n$ whenever $F$ lies inside the confidence envelope given in (2.5); hence upper bounds on overestimating $\varepsilon_n$ can be given in terms of the coverage probability of the confidence envelope. In this section, we collect a few results that are useful throughout the remainder of this paper. Readers less interested in technical ideas may prefer to skip this section and to refer back to it as needed.

A particularly easy way to analyze the confidence band given in (2.5) is through the distribution $W_n^*$ given by

$$W_n^* \stackrel{d}{=} \sup_t \left\{ \sqrt{n}\frac{|F_n(t)-F(t)|}{\sqrt{F(t)(1-F(t))}} \right\},$$

especially once we recall that the distribution of $W_n^*$ does not depend on $F$. More specifically, consider $n$ independent samples $U_i$ from a uniform distribution $U(0,1)$. The empirical distribution corresponding to these observations is then given by $V_n(t) = \frac{1}{n} \sum_{i=1}^{n} 1_{\{U_i \leq t\}}$. Set $U_n(t) = \sqrt{n}[V_n(t) - t]$, $0 < t < 1$, and write the *normalized uniform empirical process* as $W_n(t) = \frac{|U_n(t)|}{\sqrt{t(1-t)}}$. The distribution of $W_n^*$ can then be written as $W_n^* \equiv \sup_t W_n(t)$. The following well-known result [15] can be used to construct asymptotic fixed level one-sided confidence intervals for $\varepsilon_n$:

$$(3.1) \qquad \lim_{n \to \infty} \frac{W_n^*}{\sqrt{2 \log \log n}} \xrightarrow{p} 1.$$

Such an analysis underlies some of the theory in Meinshausen and Rice [14] but for the results given in our paper this approach does not suffice for reasons that we now explain.

We are interested in estimators which underestimate $\varepsilon_n$ with high probability. These estimators correspond to choosing large $a_n$ and are used to construct estimators with good mean squared error performance. Unfortunately $W_n^*$ has an extremely heavy tail [5],

$$\lim_{w \to \infty} w^2 P\{W_n^* \geq w\} = C,$$

so using $W_n^*$ to bound such tail probabilities only yields bounds on the chance that $\hat{\varepsilon}_{a_n}^*$ exceeds $\varepsilon_n$ which decrease slowly in $a_n$. Such bounds are insufficient in our analysis of the mean squared error. The reason for this is that the heavy-tailed behavior exhibited by $W_n^*$ is caused by the tails in the empirical process and in our analysis we only consider values of $t$ between $0$ and $\sqrt{2 \log n}$. Hence instead of looking at $W_n^*$ we may instead analyze the following modified version of $W_n^*$:

$$(3.2) \qquad Y_n \overset{d}{=} \max_{\{0 \leq t \leq \sqrt{2 \log n}\}} \left\{ \sqrt{n} \frac{|F_n(t) - F(t)|}{\sqrt{F(t)(1-F(t))}} \right\},$$

which can be equivalently written as $Y_n =_d \max_{\{F(0) \leq t \leq F(\sqrt{2 \log n})\}} \left\{ \frac{|U_n(t)|}{\sqrt{t(1-t)}} \right\}$.

The problem here is that $F(0)$ and $F(\sqrt{2 \log n})$ are unknown and depend on $F$, so we need a different way to estimate the tail probability of $Y_n$. We suggest two possible approaches. The first one is clean but conservative and is particularly valuable for theoretical development. The second one has a more complicated form but is sharp and allows for greater precision in the construction of confidence intervals. In the first approach, write $W_n^+$ for the distribution of $Y_n$ where $F$ corresponds to $N(0,1)$ and $F_n$ is the empirical c.d.f. formed from $n$ i.i.d. $N(0,1)$ observations. Then $W_n^+$ can be written as

$$W_n^+ \overset{d}{=} \max_{\{1/2 \leq t \leq \Phi(\sqrt{2 \log n})\}} \left\{ \frac{|U_n(t)|}{\sqrt{t(1-t)}} \right\}.$$

The following lemma shows the tail probability of any $Y_n$ associated with an $F$ is at most twice as large as that of $W_n^+$, uniformly for all Gaussian mixtures $F$ of the form $F(t) = \int \Phi(t - \mu) \, dH$ with $P\{0 \leq H \leq \sqrt{2 \log n}\} = 1$.

LEMMA 3.1. *Suppose that $Y_n$ is the distribution given in (3.2) where $F$ is a Gaussian mixture $F(t) = \int \Phi(t - \mu) \, dH$ with $P\{0 \leq H \leq \sqrt{2 \log n}\} = 1$. Then for any constant $c$, $P\{Y_n \geq c\} \leq 2 \cdot P\{W_n^+ \geq c\}$.*

The following tail bound for $W_n^+$ can be used to bound $P(\varepsilon_{a_n}^* > \varepsilon_n)$.

LEMMA 3.2. *For any constant $c_0 > 0$, for sufficiently large $n$, there is a constant $C > 0$ such that $P\{W_n^+ \geq c_0 \log^{3/2}(n)\} \leq C \cdot n^{-1.5 c_0 / \sqrt{8\pi}}$.*

It should now be clear why in our setting it is preferable to use such bounds since the corresponding tail behavior of $W_n^*$ satisfies $P\{W_n^* \geq c_0 \log^{3/2}(n)\} \asymp C \times (\log n)^{-3}$, which is not sufficient for our analysis of mean squared error given in the next section.

In the second approach, note that $F(\sqrt{2 \log n}) \leq \Phi(\sqrt{2 \log n})$, and with overwhelming probability, $F(0) \geq F_n(0) - \sqrt{c_0 \log(n)}/\sqrt{n}$. Now, for any constant $c_0 > 0$, define

$$W_n^{++} \equiv W_n^{++}(c_0) = \max_{\{(F_n(0) - \sqrt{c_0 \log(n)}/\sqrt{n}) \leq t \leq \Phi(\sqrt{2 \log n})\}} \left\{ \frac{|U_n(t)|}{\sqrt{t(1-t)}} \right\}.$$

The following lemma shows that the tail probability of $Y_n$ is almost bounded by that of $W_n^{++}$, uniformly for all one-sided mixtures even *without* the constraint that $H \leq \sqrt{2 \log n}$.

LEMMA 3.3. *Suppose that $Y_n$ is the distribution given in (3.2) where $F$ is a Gaussian mixture $F(t) = \int \Phi(t - \mu) \, dH$ with $P\{H \geq 0\} = 1$. Then for any constant $c_0 > 0$ and $c$, $P\{Y_n \geq c\} \leq P\{W_n^{++} \geq c\} + 2 n^{-c_0} \cdot (1 + o(1))$.*

This lemma is particularly useful in the construction of accurate confidence intervals where we take $c_0 = 3$ so that the difference between the two probabilities is $O(n^{-3})$. Without further notice, we refer $W_n^{++}$ to the one with $c_0 = 3$. Lemmas 3.1–3.3 are proved in [2], Sections 8.2–8.4.

3.1. *Choice of $a_n$ in later sections.* Different choices of $a_n$ lead to different estimators of $\varepsilon_n$. We shall choose $a_n$ depending on the purpose. In Section 4 the focus is on optimal rates of convergence for mean squared error. For this purpose it is convenient to choose a relatively large $a_n$ [i.e., $4 \sqrt{2\pi} \log^{3/2}(n)$]. In Section 6, where the focus is on consistency, a much smaller $a_n$ is also sufficient and might be preferred. Finally, the interest of Section 5 is on one-sided confidence intervals, and here we wish to choose an $a_n$ with level $\alpha = P\{Y_n \geq a_n\}$ being fixed. The difficulty here is that, different from the above two cases, the $a_n$ depends on the unknown $F(0)$ and $F(\sqrt{2 \log n})$. Fortunately, the level $\alpha$ is fixed and specified before hand, so one can use simulated values of $W_n^{++}$ to approximate $a_n$ without much computational complexity.

**4. Mean squared error.** In this section, we focus on choosing a member of the family of estimators constructed in Section 2.2 which has near optimal mean squared error properties. More discussion is given in Section 7 where a simulation study provides further insight into the mean squared error performance of these estimators. Our analysis begins with the bound

$$E\left(\frac{\hat{\varepsilon}_{a_n}^*}{\varepsilon_n} - 1\right)^2 \leq \left(\frac{1}{\varepsilon_n}\right)^2 P(\hat{\varepsilon}_{a_n}^* > \varepsilon_n) + E\left[\left(\frac{\hat{\varepsilon}_{a_n}^*}{\varepsilon_n} - 1\right)^2 \cdot 1_{\{\hat{\varepsilon}_{a_n}^* \leq \varepsilon_n\}}\right].$$

There is a tradeoff depending on the choice of $a_n$. As $a_n$ increases $P(\hat{\varepsilon}_{a_n}^* > \varepsilon_n)$ decreases but when $\hat{\varepsilon}_{a_n}^*$ underestimates $\varepsilon_n$ it does so by a greater amount. It is thus desirable to choose the smallest $a_n$ so that the first term is negligible and this in fact leads to an estimator with near optimal performance. It should be stressed that in the construction of the smallest such $a_n$ the precise bounds given in Lemma 3.2 are important and the tail bounds for $W_n^*$ do not suffice. In particular Lemma 3.2 shows that $a_n = 4\sqrt{2\pi} \log^{3/2}(n)$ suffices to make this first term negligible. For such a choice, the following theorem gives upper bounds on the minimax risk.

THEOREM 4.1. *Suppose* $F(t) = (1 - \varepsilon_n)\Phi(t) + \varepsilon_n \Phi(t - \mu_n)$ *with* $\varepsilon_n = n^{-\beta}$, $\mu_n = \sqrt{2r \log n}$, *where* $0 < r < 1$, $\frac{1}{2} < \beta < 1$, *and* $r > \rho^*(\beta)$ *so that* $(\beta, r)$ *falls into the interior part of the detectable region. Set* $a_n = 4\sqrt{2\pi} \log^{3/2}(n)$. *The estimator* $\hat{\varepsilon}_{a_n}^*$ *defined in* (2.8) *satisfies*

$$(4.1) \quad \begin{aligned} &E\left[\frac{\hat{\varepsilon}_{a_n}^*}{\varepsilon_n} - 1\right]^2 \\ &\leq \begin{cases} C(r, \beta)(\log n)^{5.5} n^{-1 - 2r + 2\beta}, & \text{when } \beta \geq 3r, \\ C(\beta, r)(\log n)^{5.5} n^{-1 + (\beta + r)^2/(4r)}, & \text{when } r < \beta < 3r, \\ C(r, \beta)(\log n)^4 n^{-1 + \beta}, & \text{when } \beta \leq r, \end{cases} \end{aligned}$$

*where* $C(\beta, r)$ *is a generic constant depending on* $(\beta, r)$.

Theorem 4.1 gives an upper bound for the rate of convergence of $\hat{\varepsilon}_{a_n}^*$. Although this estimator usually underestimates $\varepsilon_n$, the lower bounds for the mean squared error given below show that the performance of the estimator cannot be significantly improved.

Although the lower bounds given below are based on a two-point testing argument we should stress that they do not follow from the testing theory developed in [8]. In particular the detection boundary mentioned in Section 1 is derived by testing the simple hypothesis that $\varepsilon_n = 0$ against a particular alternative hypothesis. Here we need to study a more complicated hypothesis testing problem where both the null and alternative hypotheses correspond to Gaussian mixtures. More specifically, let $X_1, \ldots, X_n \overset{\text{i.i.d.}}{\sim} P$ and consider the following problem of testing between the two Gaussian mixtures:

$$H_0 \colon P = P_0 = (1 - \varepsilon_{0,n})N(0, 1) + \varepsilon_{0,n}N(\mu_{0,n}, 1)$$

and

$$H_1 \colon P = P_1 = (1 - \varepsilon_{1,n})N(0, 1) + \varepsilon_{1,n}N(\mu_{1,n}, 1).$$

Minimax lower bounds for estimating $\varepsilon_n$ can then be given based on carefully selected values of $\varepsilon_{0,n}$, $\varepsilon_{1,n}$, $\mu_{0,n}$ and $\mu_{1,n}$ along with good bounds on the Hellinger affinity between $n$ i.i.d. observations with distributions $P_0$ an $P_1$. As is shown in the proof of the following theorem, these bounds require somewhat delicate arguments. We should mention that our attempts using bounds on the chi-square distance, a common approach to such problems, did not yield the present results. The lower bounds are summarized as follows.

THEOREM 4.2.  Let $X_1, \ldots, X_n \overset{\text{i.i.d.}}{\sim} (1 - \varepsilon_n)N(0, 1) + \varepsilon_n N(\mu_n, 1)$. For $0 < r < 1$, $\frac{1}{2} < \beta < 1$, $a_1, a_2 > 0$ and $b_2 > b_1 > 0$, set $\Omega_n = \{(\varepsilon_n, \mu_n) \colon b_1 n^{-\beta} \leq \varepsilon_n \leq b_2 n^{-\beta}, \sqrt{2r \log n} - \frac{a_1}{\log n} \leq \mu_n \leq \sqrt{2r \log n} + \frac{a_2}{\log n}\}$.  Then

$$\inf_{\hat{\varepsilon}_n} \sup_{(\varepsilon_n, \mu_n) \in \Omega_n} E\left(\frac{\hat{\varepsilon}_n}{\varepsilon_n} - 1\right)^2 \geq \begin{cases} C(\log n)n^{-1-2r+2\beta}, & \text{when } \beta \geq 3r, \\ C(\log n)^{5/2}n^{-1+(\beta+r)^2/(4r)}, & \text{when } r < \beta < 3r, \\ Cn^{-1+\beta}, & \text{when } \beta \leq r. \end{cases}$$

A comparison between the upper bounds given in Theorem 4.1 and the lower bounds given in Theorem 4.2 shows that the procedure $\hat{\varepsilon}^*_{a_n}$ has mean squared error within a logarithmic factor of the minimax risk. Additional insight into the performance of this estimator is given in Section 6 where comparisons to an estimator introduced by Meinshausen and Rice [14] are made and in Section 7 where we report some simulations results.

**5. One-sided confidence intervals.**  In the previous section we showed how to choose $a_n$ so that the estimator $\hat{\varepsilon}^*_{a_n}$ has good mean squared error properties. In the present section we consider in more detail one-sided confidence intervals. For such intervals there are two conflicting goals. We want to maintain coverage probability over a large class of models while minimizing the amount that our estimator underestimates $\varepsilon_n$. More specifically, the goal can be formulated in terms of the following optimization problem:

$$\text{Minimize} \quad E(\varepsilon_n - \hat{\varepsilon}_n)_+ \quad \text{subject to} \quad \sup_{\mathcal{F}} P(\hat{\varepsilon}_n > \varepsilon_n) \leq \alpha,$$

where $\mathcal{F}$ is a collection of Gaussian mixtures. A similar formulation for the construction of optimal nonparametric confidence intervals is given in Cai and Low [3].

In the present section we focus on this optimization problem for the class of all two-point Gaussian mixtures showing that the estimator $\hat{\varepsilon}^*_{a_n}$ with an appropriately chosen $a_n$ provides an almost optimal lower end point for a one-sided confidence interval with a given coverage probability. Perhaps equally interesting is that this one-sided confidence interval maintains coverage probability over a much larger collection of Gaussian mixture models, namely the set of all one-sided Gaussian mixtures with $H > 0$. See also Section 6.3 where we briefly discuss how the condition $H > 0$ can be dropped.

5.1. *Coverage over one-sided Gaussian mixtures.* In this section we show how one-sided confidence intervals with a given coverage probability can be constructed for the collection of all one-sided Gaussian mixtures (1.1) with $H > 0$. Let $\mathcal{F}$ be the collection of all one-sided Gaussian mixture c.d.f.s of the form $(1-\varepsilon)\Phi(t) + \varepsilon G$ where $G(t) = \int \Phi(t-\mu)\,dH$ is the convolution of $\Phi$ and a c.d.f. $H$ supported on the positive half-line. For arbitrary constants $0 < a < b < 1$ and $0 < \tau < \tau'$, out of all c.d.f.s $F \in \mathcal{F}$ passing through points $(\tau, a)$ and $(\tau', b)$, the most "sparse" one (i.e., smallest $\varepsilon$) is a two-point Gaussian mixture $F^*(t) = (1-\varepsilon^*)\Phi(t) + \varepsilon^*\Phi(t - \mu^*)$, where $(\varepsilon^*, \mu^*)$ are chosen such that $F^*(\tau) = a$ and $F^*(\tau') = b$. That is,

(5.1)
$$\mu^*: \text{solution of } D(\mu; \tau, \tau') = \frac{\Phi(\tau) - a}{\Phi(\tau') - b} \quad \text{and}$$

$$\varepsilon^* = \frac{\Phi(\tau) - a}{\Phi(\tau) - \Phi(\tau - \mu^*)},$$

where the function $D$ is given in (2.3). The following lemma is proved in [2], Section 8.7.

LEMMA 5.1. *Fix $0 < a < b < 1$, $0 < \tau < \tau'$, and $0 < \varepsilon \leq 1$. For any $F = (1-\varepsilon)\Phi(t) + \varepsilon G \in \mathcal{F}$ such that $F(\tau) = a$ and $F(\tau') = b$, define $\varepsilon^*$ by (5.1). Then $\varepsilon^* \leq \varepsilon$.*

See Figure 1.

We now turn to the coverage probability of the grid procedure $\hat{\varepsilon}^*_{a_n}$ over the class $\mathcal{F}$. Fix an $F \in \mathcal{F}$. Then for each pair of adjacent points $(t_j, t_{j+1})$ in the grid, the above lemma shows that there is a two-point Gaussian mixture $F^*(t) = (1 - \varepsilon_j^*)\Phi(t) + \varepsilon_j^*\Phi(t - \mu_j^*)$, where $(\varepsilon_j^*, \mu_j^*)$ are chosen such that $F^*(t_j) = F(t_j)$ and $F^*(t_{j+1}) = F(t_{j+1})$. It is clear that $\varepsilon_j^*$ depends on the points $t_j$ and $t_{j+1}$, but Lemma 5.1 shows that in each case $\varepsilon_j^* \leq \varepsilon$. Now suppose that $F$ lies inside the confidence envelope defined by (2.5). In this case it follows that $\hat{\varepsilon}_{a_n}^{(j)}$ defined by (2.7) satisfies $\hat{\varepsilon}_{a_n}^{(j)} \leq \varepsilon_j^*$ and hence also $\hat{\varepsilon}_{a_n}^{(j)} \leq \varepsilon_n$. Since this holds for all $j$, it then immediately follows that $\hat{\varepsilon}_{a_n}^* \leq \varepsilon_n$ whenever $F$ lies inside the confidence interval defined by (2.5). A given level confidence interval can then be given based on the distributions of $W_n^+$ and $W_n^{++}$. This result is summarized in the following theorem.

THEOREM 5.1.   *Fix $0 < \alpha < 1$ and let $a_n$ be chosen so that $P(W_n^+ \geq a_n) \leq \alpha/2$. Then uniformly for $n$ and all one-sided Gaussian location mixtures defined in (1.2) with $P(0 < H \leq \sqrt{2\log n}) = 1$, $P\{\hat{\varepsilon}_{a_n}^* \leq \varepsilon_n\} \geq (1 - \alpha)$. Moreover, let $a_n$ be chosen so that $P(W_n^{++} \geq a_n) \leq \alpha$. Then as $n \to \infty$, uniformly for all one-sided Gaussian location mixtures defined in (1.2) with $P\{H > 0\} = 1$, $P\{\hat{\varepsilon}_{a_n}^* \leq \varepsilon_n\} \geq (1 - \alpha)(1 + o(1))$.*

5.2. *Optimality under two-point Gaussian mixture model.* In the previous section we focused on the coverage property of the one-sided confidence interval over the general class of one-sided Gaussian mixtures. In this section we return to the class of two-point Gaussian mixtures and study how "close" the lower confidence limit $\hat{\varepsilon}_n$ is to the true but unknown $\varepsilon_n$. In particular we compare the performance of our procedure with the following lower bound.

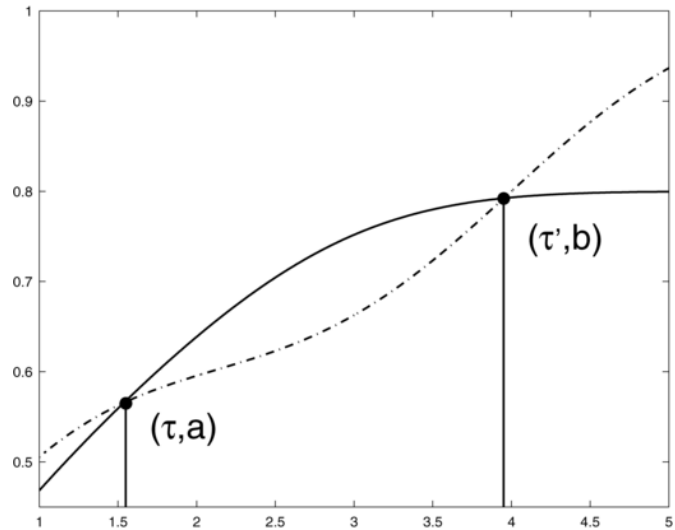

FIG. 1.   *In the c.d.f. plane, among the family of all one-sided Gaussian location mixtures which pass through two given points $(\tau, a)$ and $(\tau', b)$, the most sparse mixture is a two-point mixture (the solid curve) which bounds all other c.d.f.s from above over the whole interval $[\tau, \tau']$.*

THEOREM 5.2.    Let $X_1, \ldots, X_n \overset{\text{i.i.d.}}{\sim} (1 - \varepsilon_n)N(0, 1) + \varepsilon_n N(\mu_n, 1)$. For $0 < r < 1$, $\frac{1}{2} < \beta < 1$, $a_1, a_2 > 0$ and $b_2 > b_1 > 0$, set $\Omega_n = \{(\varepsilon_n, \mu_n) : b_1 n^{-\beta} \leq \varepsilon_n \leq b_2 n^{-\beta}, \sqrt{2r \log n} - \frac{a_1}{\log n} \leq \mu_n \leq \sqrt{2r \log n} + \frac{a_2}{\log n}\}$. For $0 < \alpha < \frac{1}{2}$, let $\hat{\varepsilon}_n$ be a $(1 - \alpha)$ level lower confidence limit for $\varepsilon_n$ over $\Omega_n$, namely, $\inf_{\Omega_n} P\{\varepsilon_n \geq \hat{\varepsilon}_n\} \geq 1 - \alpha$. Then

$$\inf_{\hat{\varepsilon}_n} \sup_{(\varepsilon_n, \mu_n) \in \Omega_n} E\left(1 - \frac{\hat{\varepsilon}_n}{\varepsilon_n}\right)_+$$

$$\geq \begin{cases} C(\log n)^{1/2} n^{-1/2 - r + \beta}, & \text{when } \beta \geq 3r, \\ C(\log n)^{5/4} n^{-1/2 + (\beta + r)^2/(8r)}, & \text{when } r < \beta < 3r, \\ C n^{-1/2 + \beta/2}, & \text{when } \beta \leq r. \end{cases}$$

Theorem 5.2 shows that even if the goal is to create an honest confidence interval over the class of two-point Gaussian mixture models the resulting estimator must underestimate the true $\varepsilon_n$ by a given amount. The following theorem shows that the estimator given in the previous section which has guaranteed coverage over the class of all one-sided Gaussian mixture models is almost optimal for two-point Gaussian mixtures.

THEOREM 5.3.    Suppose $F$ is a two-point mixture $F(t) = (1 - \varepsilon_n)\Phi(t) + \varepsilon_n \Phi(t - \mu_n)$ with $\varepsilon_n = n^{-\beta}$, $\mu_n = \sqrt{2r \log n}$, where $0 < r < 1$, $\frac{1}{2} < \beta < 1$ and $r > \rho^*(\beta)$ so $(\beta, r)$ falls into the interior part of the detectable region. Fix $0 < \alpha < 1$ and let $a_n$ be chosen so that either $P(W_n^+ \geq a_n) \leq \frac{\alpha}{2}$ or such that $P\{W_n^{++} \geq a_n\} \leq \alpha$ and for this value of $a_n$ let $\hat{\varepsilon}_{a_n}^*$ be the estimator defined in (2.8). Then there is a constant $C = C(\beta, r) > 0$ such that

$$E\left(1 - \frac{\hat{\varepsilon}_{a_n}^*}{\varepsilon_n}\right)_+$$

$$\leq \begin{cases} C \cdot \sqrt{\log \log(n)} \cdot (\log n)^{5/4} \cdot n^{-1/2 - r + \beta}, & \text{when } \beta > 3r, \\ C \cdot \sqrt{\log \log(n)} \cdot (\log n)^{5/4} \cdot n^{-1/2 + (\beta + r)^2/(8r)}, & \text{when } r < \beta \leq 3r, \\ C \cdot \sqrt{\log \log(n)} \cdot n^{-1/2 + \beta/2}, & \text{when } \beta \leq r. \end{cases}$$

**6. Discussion.**    In this section we compare and contrast the methodology developed in the present paper to the approach taken by Meinshausen and Rice [14]. The goal is to explain intuitively some of the theory developed in these two papers. Both methods have a root based on the idea of "thresholding," and how well each method works can partially be explained in terms of the concept of *most informative threshold*.

We shall start with a general comparison of the two estimators. It is useful to note that the stochastic fluctuations of these estimators are not larger in order of magnitude than the bias. It is thus instructive for a heuristic analysis to replace each of these estimators by nonrandom approximations. The approach taken in Meinshausen and Rice [14] starts with a more general mixture model which after a transformation can be written as

$$Y_i \overset{\text{i.i.d.}}{\sim} (1 - \varepsilon_n)N(0,1) + \varepsilon_n F, \qquad 1 \le i \le n,$$

where $F$ is an arbitrary distribution. In that context one-sided bounds are given for $\varepsilon_n$ which hold no matter the distribution of $F$. The lower bound can be thought of first picking an arbitrary threshold $t$, then comparing the fraction of samples $\ge t$ with the expected fraction $\ge t$ when all samples are truly from $N(0,1)$; the difference between two fractions either comes from stochastic fluctuations or from the signal, which thus naturally provides a lower bound if the stochastic fluctuations are controlled.

Using our notation, Meinshausen and Rice's lower bound can be written as $\hat{\varepsilon}^{\text{MR}}_{a_n^*} \equiv \sup_{\{-\infty < t < \infty\}} \hat{\varepsilon}^{\text{MR}}_{a_n^*}(t; F_n)$, where

$$(6.1) \qquad \hat{\varepsilon}^{\text{MR}}_{a_n^*}(t; F_n) = \left[ \frac{\Phi(t) - F_n(t) - (a_n^*/\sqrt{n}) \cdot \sqrt{\Phi(t)(1 - \Phi(t))}}{\Phi(t)} \right].$$

Here $a_n^* > 0$ is a constant which plays a similar role as $a_n$ in our estimator, and without loss of generality, we chose $1/\sqrt{t(1-t)}$ as the bounding function [14]. A useful approximation to this estimator is given by neglecting the stochastic fluctuation where we replace $F_n$ by $F$. The result is the approximation $\hat{\varepsilon}^{\text{MR}}_{a_n^*}(t; F)$,

$$
\begin{aligned}
(6.2) \qquad \hat{\varepsilon}^{\text{MR}}_{a_n^*}(t; F_n) &\approx \hat{\varepsilon}^{\text{MR}}_{a_n^*}(t; F) \\
&\equiv \left[ \frac{\Phi(t) - F(t) - (a_n^*/\sqrt{n}) \cdot \sqrt{\Phi(t)(1 - \Phi(t))}}{\Phi(t)} \right].
\end{aligned}
$$

It is instructive to compare this approximation with the following slightly modified version of our estimator where we neglect the stochastic difference by replacing $\hat{\mu}_j$ by $\mu_n$ and where we approximate $F^+$ by $F + \frac{a_n}{\sqrt{n}}\sqrt{F(1-F)}$. Then the estimator $\hat{\varepsilon}^*_{a_n}$ can be approximated by $\check{\varepsilon}^*_{a_n} \approx \sup_{\{0 \le t \le \sqrt{2\log n}\}} \check{\varepsilon}^*_{a_n}(t, F)$, where

$$(6.3) \qquad \check{\varepsilon}^*_{a_n}(t, F) = \frac{\Phi(t) - F(t) - (a_n/\sqrt{n}) \cdot \sqrt{F(t)(1 - F(t))}}{\Phi(t) - \Phi(t - \mu_n)}.$$

It is now easy to compare (6.2) with (6.3). There are three differences: (i) we use $\Phi(t) - \Phi(t - \mu_n)$ as the denominator instead of $\Phi(t)$; (ii) we use $\sqrt{F(t)(1 - F(t))}$ rather than $\sqrt{\Phi(t)(1 - \Phi(t))}$ for controlling stochastic fluctuation; (iii) we take the maximum over $(0, \sqrt{2\log n})$ instead of $(-\infty, \infty)$. In fact, only the first difference is important in the analysis of the two-point mixture model.

6.1. *Consistent estimation.* In this section we compare the approximations for the two-point mixture models starting with the Meinshausen and Rice procedure [14]. We have

$$(6.4) \quad 1 - \hat{\varepsilon}_{a_n^*}^{\mathrm{MR}}(t, F)/\varepsilon_n = \left[ \frac{\Phi(t - \mu_n)}{\Phi(t)} + a_n^* \cdot n^{\beta - 1/2} \cdot \sqrt{(1 - \Phi(t))/\Phi(t)} \right],$$

and in order for $\hat{\varepsilon}_{a_n^*}^{\mathrm{MR}}$ to be consistent, we need a $t$ such that

$$(6.5) \quad \frac{\Phi(t - \mu_n)}{\Phi(t)} \approx 0 \quad \text{and} \quad a_n^* n^{\beta - 1/2} \cdot \sqrt{(1 - \Phi(t))/\Phi(t)} \approx 0.$$

It is easy to check that both these conditions hold only if $\sqrt{2(2\beta - 1) \log n} \le t < \mu_n$ and that this is only possible when $r > 2\beta - 1$. Hence the Meinshausen and Rice procedure is only consistent on a subset of the detectable regions. Note here that consistency requires a constraint on $t$, namely that $t$ should not exceed $\mu_n$ regardless of the value of $\beta$.

A similar analysis can be provided for the approximation of our estimator. Since we use the term $\Phi(t) - \Phi(t - \mu_n)$ as the denominator in (6.3) instead of $\Phi(t)$, the above restriction on the choice of $t$ for Meinshausen and Rice's lower bound does not apply to our estimator. In fact we should always choose $t$ to be greater than $\mu_n$, not smaller; see Table 1 for the most informative $t$. This extra freedom in choosing $t$ yields the consistency over a larger range of $(\beta, r)$. In fact for the two-point Gaussian mixture model the following theorem shows that our estimator is consistent for $\varepsilon_n$ over the entire detectable region and in this sense the estimator is optimally adaptive.

THEOREM 6.1. *Let $\Omega$ be any closed set contained in the interior of the detectable region of the $\beta$-$r$ plane: $\{(\beta, r) : \rho^*(\beta) < r < 1, \frac{1}{2} < \beta < 1\}$. For any sequence of $a_n$ such that $a_n/\sqrt{2 \log \log n} \to 1$ and $P\{W_n^+ \ge a_n\}$ tends to 0, then for all $\delta > 0$,*

$$\lim_{n \to \infty} \sup_{\{(\beta, r) \in \Omega\}} P\left\{ \left| \frac{\hat{\varepsilon}_{a_n}^*}{\varepsilon_n} - 1 \right| \ge \delta \right\} = 0.$$

Figure 2 plots on the $\beta$-$r$ plane the detection boundary which separates the detectable and undetectable regions, and the classification boundary which separates classifiable and unclassifiable regions. When $(\beta, r)$ belongs to the classifiable region, it is also able to reliably tell individually which are signal and which are not. The dashed line is the separating line of consistency of the Meinshausen and Rice lower bound: above which the lower bound is consistent to $\varepsilon_n$, below which it is not; see Meinshausen and Rice [14]. The right panel of Figure 2 shows seven subregions in the detectable region as in Table 1 given in Section 6.2.

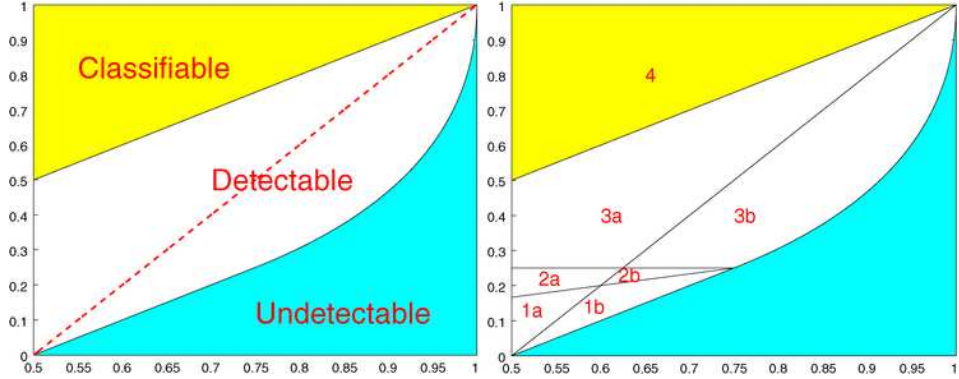

FIG. 2. Left panel: *The detection boundary and the classification boundary together with the separating line of consistency of Meinshausen and Rice (dashed line).* Right panel: *seven subregions in the detectable region as in Table 1.*

6.2. *Most informative threshold.* In this section we turn to an intuitive understanding of the mean squared error property which is driven by the value of $t$ that minimizes (6.3). More specifically, if we ignore the log-factor, the mean squared error of the estimator given by the approximation in (6.3) for a fixed $t$ satisfies

$$(1 - \hat{\varepsilon}^*_{a_n}(t, F)/\varepsilon_n)^2 = n^{2\beta - 1} \cdot \frac{F(t)(1 - F(t))}{[\Phi(t) - \Phi(t - \mu_n)]^2}.$$

Minimizing this expression over $t$ yields the optimal rate of convergence as given in Theorem 4.1. We call the minimizing value of $t$ the *most informative threshold* and these values are tabulated in Table 1. Although the mean squared error performance of the Meinshausen and Rice procedure has not been computed it appears likely that a similar phenomenon holds. In this case,

$$(1 - \hat{\varepsilon}^{\mathrm{MR}}_{a_n}(t, F)/\varepsilon_n)^2 = \left[ \frac{\Phi(t - \mu_n)}{\Phi(t)} + n^{\beta - 1/2} \cdot \sqrt{(1 - \Phi(t))/\Phi(t)} \right]^2,$$

and the value of $t$ which minimizes these expressions $\sim (2 - [2 - \frac{2\beta - 1}{r}]^{1/2})\mu_n$. Here we have assumed $r > 2\beta - 1$ as otherwise the estimator is not consistent and the most informative $t$ is not of interest; see Table 1. This shows that

$$(1 - \hat{\varepsilon}^{\mathrm{MR}}_{a_n}(t, F)/\varepsilon_n)^2 \sim n^{-(\sqrt{2r - 2\beta + 1} - \sqrt{r})^2},$$

which should give the correct convergence rate for the mean squared error. Here we have also omitted a log-factor. Since this convergence rate is always slower than the optimal rate of convergence given in Theorem 4.1, it appears at least according to this heuristic analysis that the optimal rate is never achieved by Meinshausen and Rice's estimator. One possible reason for the slow convergence rate is that in the analysis of the Meinshausen and Rice procedure the most informative $t^*$ never exceeds $\mu_n$, whereas for our procedure the most informative $t^*$ is never less than $\mu_n$. The most informative thresholds are summarized in Table 1. Note that when $r \leq 2\beta - 1$, the Meinshausen and Rice lower bound is not consistent, so the most informative threshold is not of interest (NOI). Detailed discussion on higher criticism can be found in [5].

TABLE 1
*Most informative threshold for Meinshausen and Rice's procedure and the newly proposed procedure and higher criticism of Donoho and Jin [5]. The labels of region are illustrated in the right panel in Figure 2*

| Regions in $(\beta, r)$ plane | Meinshausen and Rice | CJL | Higher criticism |
|---|---|---|---|
| 1a | $(2 - [2 - \frac{2\beta - 1}{r}]^{1/2}) \cdot \mu_n$ | $2\mu_n$ | $2\mu_n$ |
| 1b | NOI | $2\mu_n$ | $2\mu_n$ |
| 2a | $(2 - [2 - \frac{2\beta - 1}{r}]^{1/2}) \cdot \mu_n$ | $\frac{\beta + r}{2r} \cdot \mu_n$ | $2\mu_n$ |
| 2b | NOI | $\frac{\beta + r}{2r} \cdot \mu_n$ | $2\mu_n$ |
| 3a | $(2 - [2 - \frac{2\beta - 1}{r}]^{1/2}) \cdot \mu_n$ | $\frac{\beta + r}{2r} \cdot \mu_n$ | $\sqrt{2 \log n}$ |
| 3b | NOI | $\frac{\beta + r}{2r} \cdot \mu_n$ | $\sqrt{2 \log n}$ |
| 4 | $(2 - [2 - \frac{2\beta - 1}{r}]^{1/2}) \cdot \mu_n$ | $\mu_n$ | $\sqrt{2 \log n}$ |

6.3. *Extensions and generalizations.* We should stress that although the procedure presented in the present paper has better mean squared error performance than that of Meinshausen and Rice, the advantage of Meinshausen and Rice's lower bound is that it does not assume any distribution of non-null cases. In this section, we address some possible extensions of the Gaussian model which may also shed further light on the approach taken in the present paper.

Let $\{f(x; \mu) : \mu \geq 0\}$ be a family of density functions and let $X_1, \ldots, X_n$ be a random sample from a general one-sided mixture:

$$X_1, \ldots, X_n \overset{\text{i.i.d.}}{\sim} (1 - \varepsilon_n) f(x; 0) + \varepsilon_n \int f(x; \mu) \, dH(\mu), P(H > 0) = 1.$$

Two key components for the theory we developed in previous sections are: (A) among all cumulative distribution functions passing through a given pair of points $(\tau, a)$ and $(\tau', b)$, the most sparse one is a two-point mixture, and (B) the proposed estimator is optimally adaptive in estimating $\varepsilon_n$ for the family of two-point mixtures. We expect that our theory can be extended to a broad class of families where (A) and (B) hold.

We have shown in an unpublished manuscript that two conditions that suffice for (A) to hold are: (A1) the family of density functions is a strictly monotone increasing family: $f(x; \mu)/f(x)$ is increasing in $x$ for all $\mu > 0$, and (A2) $D(\mu; \tau, \tau')$ is strictly decreasing in $\mu > 0$ for any $\tau' > \tau > 0$ where $D(\mu; \tau, \tau') = \frac{F(\tau; 0) - F(\tau; \mu)}{F(\tau'; 0) - F(\tau'; \mu)}$ and $F(\cdot; \mu)$ is the c.d.f. corresponding to $f(\cdot; \mu)$.

It is interesting to note that the *two-sided* Gaussian mixture satisfies the above mentioned conditions. In fact, for $X$ from a two-sided Gaussian mixture, $|X|$ can be viewed as a one-sided mixture from the family of densities where $f(x; \mu) = \phi(x - \mu) + \phi(x + \mu) - 1$. It appears that (B) also holds in this case although we leave a more detailed analysis for future study.

**7. Simulations.** We have carried out a small-scale empirical study of the performance of our lower bound along with a comparison to Meinshausen and Rice's lower bound for sample sizes similar to those studied by Meinshausen and Rice. The purpose of the present section is only to highlight a few points that occurred consistently in our simulations. One of the points chosen in our study corresponded to $(\beta, r) = (4/7, 1/2)$. This parameter is in a region where both Meinshausen and Rice's lower bound and our lower bound are consistent. In our experiment, we simulated $n$ samples from a c.d.f. $F(t) = (1 - \varepsilon_n)\Phi(t) + \Phi(t - \mu_n)$, where $n = 10^7$, $\varepsilon_n = 10^{-4}$ and $\mu_n = \sqrt{2 \times 0.5 \times \log n} \approx 4$. The reason we chose such a large $n$ is that the signal is highly sparse. In fact, with the current $\beta$ and $n$, the number of signals is about 1000.

The experiments started by calculating $\alpha_n$-percentiles by simulation for $W_n^*$ needed for the Meinshausen and Rice procedure and for $Y_n$ for our procedure. Denote the percentiles by $a_n^*$ and $a_n$, respectively, so that $P(W_n^* \geq a_n^*) = \alpha_n$, and $P(Y_n \geq a_n) = \alpha_n$. Since $Y_n$ depends on the unknown parameter $F(0)$, we replace $Y_n$ by $W_n^{++}$ as in Lemma 3.3. The simulated data indicate that the difference between $W_n^+$ and $W_n^{++}$ is negligible and $P(W_n^+ \geq a_n) \approx \alpha_n$, so a convenient way to calculate $a_n$ is through $W_n^+$ instead of $Y_n$. We then generated 5,000 simulated values of $W_n^*$ and $W_n^+$, and calculated the values of $a_n^*$ and $a_n$ corresponding to eight chosen levels $\alpha_n = 0.5\%, 1\%, 2.5\%, 5\%, 7.5\%, 10\%, 25\%$ and $50\%$. The values are tabulated in Table 2.

Next, we laid out grid points for calculating the lower bound $\hat{\varepsilon}_{a_n}^*$. Since $2 \log n = 32.24$, we chose 33 equally-spaced grid points: $t_j = (j - 1)/\sqrt{2 \log n}$, $1 \leq j \leq 33$. We then ran 3,500 cycles of simulation.

- In each cycle we drew $n \cdot (1 - \varepsilon_n)$ samples from $N(0, 1)$ and $n \cdot \varepsilon_n$ samples from $N(\mu_n, 1)$ to approximate $n$ samples from the two-point mixture $(1 - \varepsilon_n)N(0, 1) + \varepsilon_n N(\mu_n, 1)$.
- For each $a_n$, we used the above simulated data and the grid points to calculate $\hat{\varepsilon}_{a_n}^*$.
- For each $a_n^*$, we used the simulated data to calculate $\hat{\varepsilon}_{a_n^*}^{\mathrm{MR}}$.

The results are summarized in Table 2, as well as Figure 3.

Table 2

*Comparison of our lower bound with Meinshausen and Rice's lower bound. The comparison is based on 3,500 independent cycles of simulations. In each cycle, we simulated $n = 10^7$ samples from a two-point mixture with $\varepsilon_n = 10^{-4}$ and $\mu_n = \sqrt{2 \times 0.5 \times \log n} \approx 4$. The lower bounds were calculated for each of the eight chosen $\alpha_n$-levels. The unsatisfactory performances of Meinshausen and Rice's lower bound are displayed in boldface, and are caused by its heavy-tailed behavior*

| | $\alpha_n$ | 0.005 | 0.01 | 0.025 | 0.05 | 0.075 | 0.10 | 0.25 | 0.50 |
|---|---|---|---|---|---|---|---|---|---|
| $\hat{\varepsilon}^*_{a_n}/\varepsilon_n$ | $\frac{a_n}{\sqrt{2\log n}}$ | 2.126 | 1.956 | 1.699 | 1.545 | 1.467 | 1.370 | 1.158 | 0.940 |
| | $P(\hat{\varepsilon}_n \geq \varepsilon_n)$ | 0 | 0 | 0.0014 | 0.0026 | 0.0043 | 0.0077 | 0.026 | 0.114 |
| | Maximum | 0.654 | 0.787 | 1.063 | 1.907 | 2.485 | 3.215 | 4.794 | 6.418 |
| | Mean | 0.456 | 0.477 | 0.516 | 0.544 | 0.560 | 0.583 | 0.651 | 0.776 |
| | Median | 0.450 | 0.471 | 0.508 | 0.531 | 0.546 | 0.562 | 0.608 | 0.677 |
| | Deviation | 0.045 | 0.049 | 0.062 | 0.085 | 0.1015 | 0.127 | 0.211 | 0.373 |
| | $E[\frac{\hat{\varepsilon}_n}{\varepsilon_n}-1]^2$ | 0.299 | 0.276 | 0.238 | 0.215 | 0.204 | 0.190 | 0.167 | 0.189 |
| | $E(1-\frac{\hat{\varepsilon}_n}{\varepsilon_n})_+$ | 0.545 | 0.523 | 0.485 | 0.458 | 0.442 | 0.421 | 0.364 | 0.285 |
| $\hat{\varepsilon}^{\mathrm{MR}}_{a^*_n}/\varepsilon_n$ | $\frac{a^*_n}{\sqrt{2\log n}}$ | 6.830 | 3.731 | 2.382 | 1.826 | 1.657 | 1.557 | 1.285 | 1.087 |
| | $P(\hat{\varepsilon}_n \geq \varepsilon_n)$ | 0 | 0 | 0 | 0.002 | 0.007 | 0.013 | 0.101 | 0.290 |
| | Maximum | 0.309 | 0.473 | 0.643 | 1.337 | **31.46** | **321.9** | **1113** | **1781** |
| | Mean | 0.252 | 0.374 | 0.477 | 0.5457 | 0.6017 | 0.836 | 5.644 | 26.158 |
| | Median | 0.251 | 0.373 | 0.472 | 0.537 | 0.562 | 0.579 | 0.639 | 0.739 |
| | Deviation | 0.018 | 0.027 | 0.041 | 0.065 | 0.795 | **7.765** | **43.04** | **123.2** |
| | $E[\frac{\hat{\varepsilon}_n}{\varepsilon_n}-1]^2$ | 0.560 | 0.393 | 0.276 | 0.211 | 0.791 | **60.31** | **1873** | **15814** |
| | $E(1-\frac{\hat{\varepsilon}_n}{\varepsilon_n})_+$ | 0.748 | 0.626 | 0.523 | 0.455 | 0.426 | 0.405 | 0.315 | 0.214 |

We draw attention to a number of features which showed up not only in this simulation but in our other simulations as well. First, the distribution of $\hat{\varepsilon}^*_{a_n}/\varepsilon_n$ has a relatively thin tail. Figure 3 gives histograms of $\hat{\varepsilon}^*_{a_n}/\varepsilon_n$ which show that when it does overestimate, it only overestimates by a factor of at most 5 or 6. Moreover, the chance of underestimation is in general much smaller than $\alpha_n$, sometimes even 10 times smaller, which suggests the theoretical upper bound for overestimation in Theorem 5.1 is quite conservative. For example, column 7 of Table 2 suggests for $\alpha_n = 25\%$ the empirical probability of overestimation $\approx 2.6\%$ which is roughly 10 times smaller. Finally, when it does underestimate, the amount of underestimation is reasonably small. In addition, the risk $E([\hat{\varepsilon}_{a_n}/\varepsilon_n]-1)^2$ and $E(1-[\hat{\varepsilon}_{a_n}/\varepsilon_n])_+$ are also reasonably small. We also note that Meinshausen and Rice's lower bound displays a heavy-tailed behavior; it can sometimes overestimate $\varepsilon_n$ by as much as 1,100 times.

The performance of $\hat{\varepsilon}^*_{a_n}$ is *not* very sensitive to different choice of $\alpha_n$ (or equivalently $a_n$). As $\alpha_n$ gets larger, slowly, the mean and median of $\hat{\varepsilon}^*_{a_n}$ increase, and $E([\hat{\varepsilon}_{a_n}/\varepsilon_n]-1)^2$ and $E(1-[\hat{\varepsilon}_{a_n}/\varepsilon_n])_+$ decrease, which suggest a better estimator for a larger $\alpha_n$ in a reasonable range, for example, $\alpha_n \leq 50\%$. The phenomenon can be interpreted by the thin tail property as well as that fact the chance of overestimation is slim: a larger $\alpha_n$ will not increase much of the chance of overestimation, but it will certainly boost the underestimation and in effect make the whole estimator more accurate.

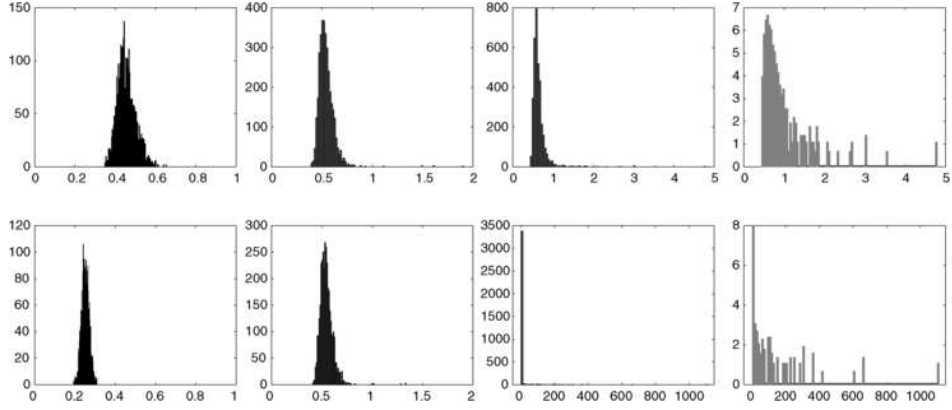

Fig. 3. *Histograms for 3,500 simulated ratios between lower bounds and the true $\varepsilon_n$. The simulations is based on $10^7$ samples from two-point mixture with $\varepsilon_n = 10^{-4}$ and $\mu_n = \sqrt{2 \times 0.5 \times \log n} \approx 4$. Top row: our lower bound. Bottom row: Meinshausen and Rice's lower bound. From left to right, lower bounds correspond to different $\alpha_n$ level: 0.005, 0.05 and 0.25. The last column is the log-histogram of the third column.*

We now turn to Meinshausen and Rice's lower bound. $\hat{\varepsilon}_{a_n^*}^{\mathrm{MR}}$ also provides an honest lower bound, and $P(\hat{\varepsilon}_{a_n^*}^{\mathrm{MR}} \geq \varepsilon_n)$ is typically much smaller than $\alpha_n$. However, for relatively larger $\alpha_n$, empirical study shows that $\hat{\varepsilon}_{a_n^*}^{\mathrm{MR}}$ is not an entirely satisfactory lower bound as the variance of $\hat{\varepsilon}_{a_n^*}^{\mathrm{MR}}$ is relatively large. For example, when $\alpha_n \geq 0.1$, $E(\frac{\hat{\varepsilon}_{a_n^*}^{\mathrm{MR}}}{\varepsilon_n} - 1)^2$ can be as large as a few hundred or a few thousand; see the cells in boldface in the table. Even for smaller $\alpha_n$, $\hat{\varepsilon}_{a_n^*}^{\mathrm{MR}}$ is slightly worse than $\hat{\varepsilon}_{a_n}^*$ if we compare the mean, median, $E([\hat{\varepsilon}_{a_n^*}^{\mathrm{MR}}/\varepsilon_n] - 1)^2$ and risks, and so on, which suggests $\hat{\varepsilon}_{a_n^*}^{\mathrm{MR}}$ is not as accurate as $\hat{\varepsilon}_{a_n}^*$.

The large variance of $\hat{\varepsilon}_{a_n^*}^{\mathrm{MR}}$ is caused by its heavy-tailed behavior. We have plotted the histograms of $\hat{\varepsilon}_{a_n^*}^{\mathrm{MR}}/\varepsilon_n$. In some circumstances, $\hat{\varepsilon}_{a_n^*}^{\mathrm{MR}}$ can overestimate $\varepsilon_n$ by a factor of several hundred or even larger, and a larger-scale study shows that this phenomenon will not disappear just by taking a smaller $\alpha_n$.

Naturally, one wonders what causes such heavy-tailed behavior and how to modify $\hat{\varepsilon}_{a_n^*}^{\mathrm{MR}}$ such that it preserves the good property of $\hat{\varepsilon}_{a_n^*}^{\mathrm{MR}}$ and with a relatively thin tail. Recall that ([14])

$$(7.1) \qquad \hat{\varepsilon}_{a_n^*}^{\mathrm{MR}} = \sup_{0 < t < 1} \left\{ \frac{F_n(t) - t - (a_n^*/\sqrt{n}) \cdot \sqrt{t(1-t)}}{1 - t} \right\};$$

the heavy-tailed behavior of $\hat{\varepsilon}_{a_n^*}^{\mathrm{MR}}$ is mainly caused by the denominator term $(1-t)$, which can become extremely small as $t$ gets closer to 1. We recommend dropping the term in the denominator and using the following as a lower bound:

$$\hat{\varepsilon}_{a_n^*}^+ = \sup_{0 < t < 1}[F_n(t) - t - (a_n^*/\sqrt{n}) \cdot \sqrt{t(1-t)}].$$

Clearly this is still a lower bound which is a little bit more conservative than $\hat{\varepsilon}_{a_n^*}^{\mathrm{MR}}$. However whenever the maximum in (7.1) is reached at $t \approx 0$, the difference between $\hat{\varepsilon}_{a_n^*}^{\mathrm{MR}}$ and $\hat{\varepsilon}_{a_n^*}^+$ is small. The advantage of this procedure is that it has a thin tail.

## 8. Proofs.

8.1. *Proof of Theorem* 4.1. Before going into technical details, we briefly explain the main ideas behind the proof. First note that there are two major contributions to the risk: one part due to overestimating $\varepsilon_n$ and the other part due to underestimating $\varepsilon_n$. By selecting $a_n$ as large as $4\sqrt{2\pi}\log^{3/2}(n)$, the probability of overestimating is so small that the first part is negligible. It is thus sufficient to limit our attention to the event where the estimator underestimates $\varepsilon_n$. Now recall that the estimator $\varepsilon_{a_n}^*$ is the maximum of a collection of *individual* estimators $\varepsilon_{a_n}^{(j)}$, each of which is based on a pair of adjacent grid points $t_j$ and $t_{j+1}$. Comparing $\varepsilon_{a_n}^*$ with $\varepsilon_{a_n}^{(j)}$, it is clear that the component of the risk due to $\varepsilon_{a_n}^*$ underestimating $\varepsilon$ will not exceed that of any $\varepsilon_{a_n}^{(j)}$; hence we can choose any such estimator to give us an upper bound for this component of the risk.

In detail, let $t_n^* = \sqrt{2q\log n}$ with

$$(8.1) \qquad q = \begin{cases} 4r, & \beta \geq 3r, \\ (\beta+r)^2/(4r), & r < \beta < 3r, \\ r, & \beta \leq r. \end{cases}$$

The particular $j = j_0$ we would like to choose is the one which satisfies $t_{j_0} \leq t_n^* < t_{j_0+1}$. To elaborate the above observations, we denote the event $\{F_{a_n}^-(t) \leq F(t) \leq F_{a_n}^+(t), \forall 0 \leq t \leq \sqrt{2\log n}\}$ by $A^{a_n}$. First, note that for $a_n = 4\sqrt{2\pi}\log^{3/2}(n)$, Lemma 3.2 implies that $P((A^{a_n})^c) = O(1/n^3)$. It then follows that in the bound for the risk given by $E(\frac{\hat{\varepsilon}_{a_n}}{\varepsilon_n} - 1)^2 \leq (\frac{1}{\varepsilon_n})^2 P((A^{a_n})^c) + E([\hat{\varepsilon}_{a_n}/\varepsilon_n - 1]^2 \cdot 1_{\{A^{a_n}\}})$, the first term is negligible. Second, note that $\hat{\varepsilon}_{a_n}^{(j_0)} \leq \hat{\varepsilon}_{a_n}^* \leq \varepsilon_n$ over $A^{a_n}$, so

$$(8.2) \qquad E([\hat{\varepsilon}_{a_n}^*/\varepsilon_n - 1]^2 \cdot 1_{\{A^{a_n}\}}) \leq E([\hat{\varepsilon}_{a_n}^{(j_0)}/\varepsilon_n - 1]^2 \cdot 1_{\{A^{a_n}\}}).$$

Finally, the key inequality we need to show is

(8.3)
$$E([\hat{\varepsilon}_{a_n}^{(j_0)}/\varepsilon_n - 1]^2 \cdot 1_{\{A^{a_n}\}})$$
$$\leq \begin{cases} C \log(n)(a_n^2/n) \dfrac{F(t_n^*)(1 - F(t_n^*))}{(\Phi(t_n^*) - F(t_n^*))^2}, & \beta > r, \\[3ex] C(a_n^2/n) \dfrac{F(t_n^*)(1 - F(t_n^*))}{(\Phi(t_n^*) - F(t_n^*))^2}, & \beta \leq r. \end{cases}$$

In fact, Theorem 4.1 follows directly by combining (8.2)–(8.3) with the following lemma in which we calculate $[F(t_n^*)(1 - F(t_n^*))]/[\Phi(t_n^*) - F(t_n^*)]^2$.

LEMMA 8.1. *Suppose* $F(\cdot) = (1 - \varepsilon_n)\Phi(\cdot) + \varepsilon_n\Phi(\cdot - \mu_n)$ *with* $\varepsilon_n = n^{-\beta}$, $\mu_n = \sqrt{2r \log n}$, *where* $1/2 < \beta < 1$, *and* $r > \rho^*(\beta)$ *so* $(\beta, r)$ *falls above the detection boundary. With* $t_n^*$ *defined in* (8.1),

$$\frac{F(t_n^*)(1 - F(t_n^*))}{[\Phi(t_n^*) - F(t_n^*)]^2}$$
$$= \begin{cases} \sqrt{\pi r \log n} \cdot n^{-2r+2\beta} \cdot (1 + o(1)), & \beta > 3r, \\[2ex] \dfrac{\beta(\beta - r)}{\beta + r}\sqrt{4\pi r \log n} \cdot n^{(\beta+r)^2/(4r)} \cdot (1 + o(1)), & r < \beta \leq 3r, \\[2ex] 2 \cdot n^\beta \cdot (1 + o(1)), & \beta \leq r. \end{cases}$$

*Moreover, for any* $|t - t_n^*| \leq c/\sqrt{\log n}$, *there is a constant* $C = C(r, \beta; c) > 0$ *such that* $F(t)(1 - F(t))/(\Phi(t) - F(t))^2 \leq C \cdot F(t_n^*)(1 - F(t_n^*))/(\Phi(t_n^*) - F(t_n^*))^2$.

Using $1 - \Phi(x) \sim \phi(x)/x$ for large $x$, the proof for Lemma 8.1 follows from basic calculus and is thus omitted.

The proof of (8.3) needs careful analysis on $|F_{a_n}^\pm - F|$ and $|\hat{\mu}_{a_n}^{(j_0)} - \mu_n|$. The following lemmas are proved in [2], Sections 8.5.1 and 8.5.2, respectively.

LEMMA 8.2. *For fixed* $0 < q < 1$, $a_n = O(\log^{3/2} n)$ *and* $t = t_n = \sqrt{2q \log n} + O(1/\sqrt{\log(n)})$, *we have that* $|F_{a_n}^\pm(t) - F(t)| \leq (a_n/\sqrt{n}) \cdot \sqrt{F(t)(1 - F(t))} \cdot (1 + o(1))$ *over the event* $A^{a_n}$.

LEMMA 8.3. *Suppose* $F(\cdot) = (1 - \varepsilon_n)\Phi(\cdot) + \varepsilon_n\Phi(\cdot - \mu_n)$ *with* $\varepsilon_n = n^{-\beta}$, $\mu_n = \sqrt{2r \log n}$, *where* $1/2 < \beta < 1$, *and* $r > \rho^*(\beta)$ *so* $(\beta, r)$ *falls above the detection boundary. Then there is a constant* $C > 0$ *such that over event* $A^{a_n}$, $\hat{\mu}_{a_n}^{(j_0)} \geq \mu_n$ *and for sufficiently large* $n$, $|\hat{\mu}_{a_n}^{(j_0)} - \mu_n| \leq C \cdot (a_n/\sqrt{n}) \cdot \sqrt{F(t_{j_0})(1 - F(t_{j_0}))}/[\Phi(t_{j_0}) - F(t_{j_0})]$. *As a result,* $E[(\hat{\mu}_{a_n}^{(j_0)} - \mu_n) \cdot 1_{\{A^{a_n}\}}]^2 \leq C \cdot (a_n^2/n) \cdot F(t_{j_0})(1 - F(t_{j_0}))/[\Phi(t_{j_0}) - F(t_{j_0})]^2$.

We now proceed to prove (8.3). For short, denote $A = A^{a_n}$, $\tau = t_{j_0}$, $\mu = \mu_n$, $\hat{\mu} = \hat{\mu}_{a_n}^{(j_0)}$, $\varepsilon = \varepsilon_n$, $\hat{\varepsilon} = \hat{\varepsilon}_{a_n}^{(j_0)}$ and $F^{\pm} = F_{a_n}^{\pm}$. By basic algebra, we can rewrite $\hat{\varepsilon}/\varepsilon - 1 = \frac{\Phi(\tau) - \Phi(\tau-\mu)}{\Phi(\tau) - \Phi(\tau-\hat{\mu})} \cdot [\frac{F(\tau) - F^+(\tau)}{\Phi(\tau) - F(\tau)} - \frac{\Phi(\tau-\mu) - \Phi(\tau-\hat{\mu})}{\Phi(\tau) - \Phi(\tau-\mu)}]$. But by Lemma 8.3, $\hat{\mu} \geq \mu$ over $A$ so the first term $\leq 1$. We then have

$$(8.4) \quad (\hat{\varepsilon}/\varepsilon - 1)^2 \leq 2\left[ \left( \frac{F(\tau) - F^+(\tau)}{\Phi(\tau) - F(\tau)} \right)^2 + \left( \frac{\Phi(\tau - \hat{\mu}) - \Phi(\tau - \mu)}{\Phi(\tau) - \Phi(\tau - \mu)} \right)^2 \right].$$

Now, first, by Lemma 8.2,

$$(8.5) \qquad E\left[ \left( \frac{F(\tau) - F^+(\tau)}{\Phi(\tau) - F(\tau)} \right)^2 \cdot 1_A \right] \sim (a_n^2/n) \cdot \frac{F(\tau)(1 - F(\tau))}{(\Phi(\tau) - F(\tau))^2},$$

and second, observe that $\frac{|\Phi(\tau-\hat{\mu}) - \Phi(\tau-\mu)|}{\Phi(\tau) - \Phi(\tau-\mu)} \sim \frac{\phi(\tau-\mu)}{\Phi(\tau) - \Phi(\tau-\mu)} \cdot |\hat{\mu} - \mu|$, where $\frac{\phi(\tau-\mu)}{\Phi(\tau) - \Phi(\tau-\mu)} = O(\tau - \mu)$ when $\beta > r$ and $= O(1)$ when $\beta \leq r$. So by Lemma 8.3,

$$
\begin{aligned}
(8.6) \quad & E\left( \left[ \frac{\Phi(\tau - \hat{\mu}_n) - \Phi(\tau - \mu_n)}{\Phi(\tau) - \Phi(\tau - \mu)} \right]^2 1_A \right) \\
& \leq \begin{cases} C \log(n)(a_n^2/n) \dfrac{F(\tau)(1 - F(\tau))}{(\Phi(\tau) - F(\tau))^2}, & \beta > r, \\[2ex] C(a_n^2/n) \dfrac{F(\tau)(1 - F(\tau))}{(\Phi(\tau) - F(\tau))^2}, & \beta \leq r; \end{cases}
\end{aligned}
$$

inserting (8.5)–(8.6) into (8.4) gives (8.3) and completes the proof of the theorem.

### 8.2. *Proof of Theorem 4.2.*

The basis strategy underlying the proof of Theorem 4.2 is to calculate the Hellinger affinity between pairs of carefully chosen probability measures since, as Le Cam and Yang [12] have shown, corresponding bounds for the minimax mean squared error easily follow. More specifically, let $Q_{\theta_1}$ and $Q_{\theta_2}$ be a pair of probability measures. The Hellinger affinity is defined by $A(Q_{\theta_1}, Q_{\theta_2}) = \int \sqrt{dQ_{\theta_1} \, dQ_{\theta_2}}$ and the minimax risk is bounded as

$$(8.7) \qquad \inf_{\hat{\theta}} \sup_{\theta \in \{\theta_1, \theta_2\}} E(\hat{\theta} - \theta)^2 \geq \tfrac{1}{16}(\theta_2 - \theta_1)^2 A^4(Q_{\theta_1}, Q_{\theta_2}).$$

The actual implementation of this general strategy in the proof of Theorem 4.2 requires great care in the choice of the two probability measures and involves somewhat delicate calculations of the affinity between these measures. Let $X_1, \ldots, X_n \overset{\text{i.i.d.}}{\sim} P$. Let $P_0 = (1 - \varepsilon_{0,n})N(0, 1) + \varepsilon_{0,n}N(\mu_{0,n}, 1)$ and $P_1 = (1 - \varepsilon_{1,n})N(0, 1) + \varepsilon_{1,n}N(\mu_{1,n}, 1)$. We shall write $\varepsilon_i$ for $\varepsilon_{i,n}$ and $\mu_i$ for $\mu_{i,n}$ for $i = 0, 1$, and calibrate by $\varepsilon_0 = n^{-\beta}$, $\varepsilon_1 = n^{-\beta} + (\log n)^{\rho} n^{-\tau}$ with $\tau \geq \beta$ and $\frac{1}{2} < \beta \leq 1$, $\mu_{0,n} = \sqrt{2r \log n}$ for some $r > 0$, and $\mu_{1,n} = \sqrt{2r \log n} - \delta_n$ where $\delta_n$ is "small" and will be specified later.

Denote by $P_{i,n}$ the joint distribution of $X_1, \ldots, X_n$ under $H_i$ for $i = 0, 1$. Set $\lambda_n = \frac{\beta + r}{2\sqrt{r}}\sqrt{2\log n}$ and $\Delta(x) \equiv \varepsilon_0(e^{\mu_0 x - \mu_0^2/2} - 1) + \varepsilon_1(e^{\mu_1 x - \mu_1^2/2} - 1) + \varepsilon_0\varepsilon_1(e^{\mu_0 x - \mu_0^2/2} - 1)(e^{\mu_1 x - \mu_1^2/2} - 1)$. Then simple calculations show that the Hellinger affinity between $P_0$ and $P_1$ satisfies

$$A(P_0, P_1) = \int \sqrt{dP_0 \, dP_1} = \int_{-\infty}^{\infty} \sqrt{1 + \Delta(x)}\phi(x)\, dx$$

$$= \left\{\int_{-\infty}^{\lambda_n} + \int_{\lambda_n}^{\infty}\right\}\sqrt{1 + \Delta(x)}\phi(x)\, dx.$$

It then follows from the inequalities $\sqrt{1 + \Delta} \geq 1 + \frac{1}{2}\Delta - \frac{1}{8}\Delta^2 + \frac{1}{16}\Delta^3 - \frac{5}{128}\Delta^4$ and $1 + \Delta(x) \geq [1 + (\varepsilon_0\varepsilon_1)^{1/2}(e^{\mu_0 x - \mu_0^2/2} - 1)^{1/2}(e^{\mu_1 x - \mu_1^2/2} - 1)^{1/2}]^2$ and some algebra that

$$A(P_0, P_1) \geq 1 - \tfrac{1}{2}\Delta_1 - \tfrac{1}{8}\Delta_2 + o(n^{-1}),$$

where

$$\Delta_1 = \varepsilon_0\tilde{\Phi}(\lambda_n - \mu_0)\left\{\left(1 - \left(\frac{\varepsilon_1}{\varepsilon_0}\right)^{1/2}\left(\frac{\tilde{\Phi}(\lambda_n - \mu_1)}{\tilde{\Phi}(\lambda_n - \mu_0)}\right)^{1/2}\right)^2\right.$$

$$+ 2\left(\frac{\varepsilon_1}{\varepsilon_0}\right)^{1/2}\left(\frac{\tilde{\Phi}(\lambda_n - \mu_1)}{\tilde{\Phi}(\lambda_n - \mu_0)}\right)^{1/2}$$

$$\left.\times \left(1 - e^{-(1/8)(\mu_0 - \mu_1)^2}\frac{\tilde{\Phi}(\lambda_n - (\mu_0 + \mu_1)/2)}{\tilde{\Phi}^{1/2}(\lambda_n - \mu_0)\tilde{\Phi}^{1/2}(\lambda_n - \mu_1)}\right)\right\},$$

$$\Delta_2 = \varepsilon_0^2 e^{\mu_0^2}\Phi(\lambda_n - 2\mu_0)$$

$$\times \left\{\left(1 - \frac{\varepsilon_1}{\varepsilon_0}e^{(1/2)(\mu_1^2 - \mu_0^2)}\left(\frac{\Phi(\lambda_n - 2\mu_1)}{\Phi(\lambda_n - 2\mu_0)}\right)^{1/2}\right)^2\right.$$

$$+ 2\frac{\varepsilon_1}{\varepsilon_0}e^{(1/2)(\mu_1^2 - \mu_0^2)}\left(\frac{\Phi(\lambda_n - 2\mu_1)}{\Phi(\lambda_n - 2\mu_0)}\right)^{1/2}$$

$$\left.\times \left(1 - e^{-(1/2)(\mu_1 - \mu_0)^2}\frac{\Phi(\lambda_n - (\mu_0 + \mu_1))}{\Phi^{1/2}(\lambda_n - 2\mu_0)\Phi^{1/2}(\lambda_n - 2\mu_1)}\right)\right\}.$$

CASE 1.   $\beta \geq 3r$. In this case set $\tau = \frac{1}{2} + r$, $\rho = \frac{1}{2}$ and $\delta_n = (2r)^{-1/2} \times n^{-\tau + \beta}$. With these choices, direct calculations show that $\Delta_2 \gg \Delta_1$ and it suffices to focus attention on $\Delta_2$ in this case. We shall only consider the case $\beta > 3r$ as the case $\beta = 3r$ is similar. When $\beta > 3r$, $\frac{\beta + r}{2\sqrt{r}} > 2\sqrt{r}$ and $\lambda_n > 2\mu_i$, $i = 0, 1$, for sufficiently large $n$. Hence $\Delta_2 = \varepsilon_0^2 e^{\mu_0^2}\{(1 - (1 + (\log n)^{\rho}n^{-\tau + \beta})(1 - \mu_0\delta_n + \frac{1}{2}\delta_n^2))^2 + 2[1 - (1 - \frac{1}{2}\delta_n^2)]\}(1 + o(1)) = \frac{1}{2r}n^{-1}(1 +$ $o(1)$). Thus $A(P_0, P_1) \geq 1 - \frac{1}{2}\Delta_1 - \frac{1}{8}\Delta_2 + o(n^{-1}) = 1 - \frac{1}{16r}n^{-1}(1 + o(1))$ and consequently $A(P_{0,n}, P_{1,n}) = A^n(P_0, P_1) \geq (1 - \frac{1}{16r}n^{-1}(1 + o(1)))^n \rightarrow e^{-1/(16r)} > 0$. It then follows that the minimax lower bound for estimation under the mean squared error satisfies $\inf_{\hat{\varepsilon}_n} \sup_{(\varepsilon_n, \mu_n) \in \Omega_n} E(\hat{\varepsilon}_n - \varepsilon_n)^2 \geq C(\varepsilon_{0,n} - \varepsilon_{1,n})^2 = C(\log n)n^{-1-2r}$ for some constant $C > 0$. Hence $\inf_{\hat{\varepsilon}_n} \sup_{(\varepsilon_n, \mu_n) \in \Omega_n} E(\frac{\hat{\varepsilon}_n}{\varepsilon_n} - 1)^2 \geq C(\log n)n^{-1-2r+2\beta}$.

CASE 2. $r < \beta < 3r$. In this case set $\tau = \frac{1}{2} + \beta - \frac{(\beta+r)^2}{8r}$, $\rho = \frac{5}{4}$ and $\delta_n = \frac{(\log n)^\rho n^{-\tau+\beta}}{\lambda_n - \mu_0} = \frac{\sqrt{2r}}{\beta - r}(\log n)^{3/4}n^{-\tau+\beta}$. Note that for sufficiently large $n$, $\mu_i < \lambda_n < 2\mu_i$ for $i = 0, 1$. In this case $\Delta_1$ and $\Delta_2$ are balanced. It then follows from the standard approximation to the Gaussian tail probability, $\tilde{\Phi}(x) = \frac{1}{\sqrt{2\pi}x}e^{-(1/2)x^2}(1 + o(1))$ as $x \rightarrow \infty$, that

$$\Delta_1 = \frac{1}{4}\varepsilon_0\tilde{\Phi}(\lambda_n - \mu_0)$$
$$\times \left\{ \left( [(\log n)^\rho n^{-\tau+\beta} - (\lambda_n - \mu_0)\delta_n] - \frac{\delta_n}{\lambda_n - \mu_0} \right)^2 + \frac{\delta_n^2}{(\lambda_n - \mu_0)^2} \right\}$$
$$\times (1 + o(1))$$
$$= \frac{1}{2}\varepsilon_0\tilde{\Phi}(\lambda_n - \mu_0)\frac{\delta_n^2}{(\lambda_n - \mu_0)^2}(1 + o(1)) = \frac{2r^{5/2}}{\sqrt{\pi}(\beta - r)^5}n^{-1}(1 + o(1))$$

and

$$\Delta_2 = \varepsilon_0^2 e^{\mu_0^2}\Phi(\lambda_n - 2\mu_0)\frac{2\delta_n^2}{(2\mu_0 - \lambda_n)^2}(1 + o(1))$$
$$= \frac{8r^{5/2}}{\sqrt{\pi}(\beta - r)(3r - \beta)^3}n^{-1}(1 + o(1)).$$

Hence $A(P_0, P_1) \geq 1 - \frac{1}{2}\Delta_1 - \frac{1}{8}\Delta_2 + o(n^{-1}) = 1 - cn^{-1}(1 + o(1))$, where $c = \frac{r^{5/2}}{\sqrt{\pi}(\beta-r)^5} + \frac{r^{5/2}}{\sqrt{\pi}(\beta-r)(3r-\beta)^3}$. Therefore $A(P_{0,n}, P_{1,n}) = A^n(P_0, P_1) \geq (1 - cn^{-1})^n \rightarrow e^{-c} > 0$ and consequently $\inf_{\hat{\varepsilon}_n} \sup_{(\varepsilon_n, \mu_n) \in \Omega_n} E(\frac{\hat{\varepsilon}_n}{\varepsilon_n} - 1)^2 \geq C(\varepsilon_{0,n} - \varepsilon_{1,n})^2 \geq C(\log n)^{5/2}n^{-1-2\beta+(\beta+r)^2/(4r)}$.

CASE 3. $\beta \leq r$. In this case set $\tau = \frac{1}{2} + \frac{1}{2}\beta$, $\rho = 0$ and $\delta_n = 0$. With these choices $\mu_0 = \mu_1$ and this case is simpler than the other two cases. It is easy to verify that $\Delta_1 \gg \Delta_2$ and $A(P_{0,n}, P_{1,n}) = A^n(P_0, P_1) \geq (1 - cn^{-1})^n \rightarrow e^{-c} > 0$ and once again it follows from (8.7) that $\inf_{\hat{\varepsilon}_n} \sup_{(\varepsilon_n, \mu_n) \in \Omega_n} E(\frac{\hat{\varepsilon}_n}{\varepsilon_n} - 1)^2 \geq Cn^{-1+\beta}$.

8.3. *Proof of Theorem* 5.1. Consider the event $A_n^{a_n} \equiv \{F_{a_n}^-(t) \leq F(t) \leq F_{a_n}^+(t) : \forall\, 0 \leq t \leq \sqrt{2\log n}\}$. For the first claim, on one hand, the above argument shows that $\varepsilon_{a_n}^* \leq \varepsilon_n$ over $A_n^{a_n}$. On the other hand, it follows directly from the definition of $F_{a_n}^{\pm}$ that $Y_n \leq a_n$ over $A_n^{a_n}$, so by Lemma 3.2, $P((A_n^{a_n})^c) \leq P(Y_n \geq a_n) \leq 2P(W^+ \geq a_n) \leq \alpha$. Combining these, the first claim follows from Lemma 5.1 and the argument right below it in Section 5. The second claim follows similarly by using Lemma 3.3.

8.4. *Proof of Theorem* 5.2. We give only a sketch of the proof of Theorem 5.2 since the details in terms of calculating the Hellinger affinity are similar to the proof of Theorem 4.2. Without loss of generality assume $b_1 \leq 1 < b_2$. Set

$$
\begin{cases}
\tau = \dfrac{1}{2} + r, & \rho = \dfrac{1}{2}, \delta_n = (2r)^{-1/2} n^{-\tau+\beta} \\
& \text{when } \beta \geq 3r, \\[2mm]
\tau = \dfrac{1}{2} + \beta - \dfrac{(\beta+r)^2}{8r}, & \rho = \dfrac{5}{4}, \delta_n = \dfrac{\sqrt{2r}}{\beta-r}(\log n)^{3/4} n^{-\tau+\beta} \\
& \text{when } r < \beta < 3r, \\[2mm]
\tau = \dfrac{1}{2} + \dfrac{1}{2}\beta, & \rho = 0, \delta_n = 0 \\
& \text{when } \beta \leq r.
\end{cases}
$$

For $\frac{1}{2} < \beta < 1$ and $0 < r < 1$, set $(\varepsilon_{0,n}, \mu_{0,n}) = (n^{-\beta}, \sqrt{2r\log n})$ and $(\varepsilon_{1,n}, \mu_{1,n}) = (\varepsilon_{0,n} + c_*(\log n)^\rho n^{-\tau}, \mu_{0,n} - \delta_n)$. It is clear that $(\varepsilon_{0,n}, \mu_{0,n})$ and $(\varepsilon_{1,n}, \mu_{1,n})$ are both in $\Omega_n$. Calculations as given in the proof of Theorem 4.2 then yield lower bounds on the Hellinger affinity which in turn give upper bounds on the $L_1$ distance between $P_{0,n}$ and $P_{1,n}$. These bounds show that for any given $0 < \gamma < \frac{1}{2}$ one can choose a constant $c_* > 0$ such that the $L_1$ distance between the distributions satisfies $L_1(P_{0,n}, P_{1,n}) \leq 2\gamma$. Since $\hat{\varepsilon}_n$ is a $(1-\alpha)$ level lower confidence limit over $\Omega_n$, $P_{0,n}(\hat{\varepsilon}_n \leq \varepsilon_{0,n}) \geq 1 - \alpha$. It then follows that $P_{1,n}(\hat{\varepsilon}_n \leq \varepsilon_{0,n}) \geq 1 - \alpha - \gamma$ and hence $E_{1,n}(\varepsilon_{1,n} - \hat{\varepsilon}_n)_+ \geq (1 - \alpha - \gamma)(\varepsilon_{1,n} - \varepsilon_{0,n}) = (1 - \alpha - \gamma)c_*(\log n)^\rho n^{-\tau}$.

8.5. *Proof of Theorem* 5.3. We will only show the first claim since the proof of the second claim is similar. Let $A_n$ be the event that $\sqrt{n}|F_n(t) - F(t)|/\sqrt{F(t)(1-F(t))} \leq 4\sqrt{2\pi}\log^{3/2}(n)$ for all $0 \leq t \leq \sqrt{2\log n}$; by Lemma 3.2 the risk over $A_n^c$ is negligible. Adapting the notation of the proof of Theorem 4.1, the key for the proof is that, similarly to the proof of Lemma 4.1, especially (8.3) and Lemma 8.1, the following is true for a wide range of $a_n$, for example, $O(\sqrt{\log\log n}) \leq a_n \leq 4\sqrt{2\pi}\log^{3/2}(n)$:

$$
(8.8) \qquad E\left(\left[1 - \frac{\hat{\varepsilon}_{a_n}^{(j_0)}}{\varepsilon_n}\right]^2 \cdot 1_{\{A_n\}}\right)
$$

$$\leq \begin{cases} Ca_n^2(\log n)^{2.5}n^{-1-2r+2\beta}, & \text{when } \beta \geq 3r, \\ Ca_n^2(\log n)^{2.5}n^{-1+(\beta+r)^2/(4r)}, & \text{when } r < \beta < 3r, \\ Ca_n^2(\log n)n^{-1+\beta}, & \text{when } \beta \leq r. \end{cases}$$

Using Hölder's inequality, and noting that $[1 - \hat{\varepsilon}_{a_n}^*/\varepsilon_n]_+ \leq [1 - \hat{\varepsilon}_{a_n}^{(j_0)}/\varepsilon_n]_+$, all we need to show is that $a_n \leq O(\sqrt{2\log\log n})$. Choose $a_n^*$ such that $P(W_n^* \geq a_n^*) = \alpha/2$, compare it with $P(W_n^+ \geq a_n) = \alpha/2$, as $W_n^+ \leq W_n^*$, so $a_n \leq a_n^*$. It is well known that $a_n^* \sim \sqrt{2\log\log n}$ for any fixed $0 < \alpha < 1$ (see, e.g., [15], page 600), so the claim follows directly.

8.6. *Proof of Theorem* 6.1. By Lemma 3.2, uniformly, the probability of over-estimation will not exceed $P\{Y_n \geq a_n\} \leq 2P\{W_n^+ \geq a_n\}$, which tends to 0 by the choice of $a_n$. So it is sufficient to show that $(1 - \hat{\varepsilon}_{a_n}^*/\varepsilon_n)_+$ tends to 0 in probability uniformly for all $(\beta, r) \in \Omega$.

Note that Theorem 5.3 still holds if we replace the sequence $a_n$ there by the current one. Moreover, the inequality can be further strengthened into a constant $C(\Omega) > 0$ such that for sufficiently large $n$

$$(8.9) \quad \begin{aligned} &E\left[\left(1 - \frac{\hat{\varepsilon}_{a_n}^*}{\varepsilon_n}\right)_+\right] \\ &\leq \begin{cases} C(\Omega)\sqrt{\log\log n} \cdot (\log n)^{5/4} \cdot n^{-[1/2+r-\beta]}, \\ \qquad \text{when } \beta \geq 3r, \\ C(\Omega)\sqrt{\log\log n} \cdot (\log n)^{5/4} \cdot n^{-[1/2-(\beta+r)^2/(8r)]}, \\ \qquad \text{when } r < \beta < 3r, \\ C(\Omega)\sqrt{\log\log n} \cdot n^{-[1/2-\beta/2]}, \\ \qquad \text{when } \beta \leq r. \end{cases} \end{aligned}$$

At the same time, note that the exponents are bounded away from 0:

$$(8.10) \quad d(\Omega) \equiv \min_\Omega\left\{\frac{1}{2} + r - \beta, \frac{1}{2} - \frac{(\beta+r)^2}{8r}, \frac{1-\beta}{2}\right\} > 0.$$

Combining (8.9) and (8.10) yields that $E[(1 - \hat{\varepsilon}_{a_n}^*/\varepsilon_n)_+] \leq C(\Omega) \cdot \sqrt{\log\log n} \times \log^{1.25}(n) \cdot n^{-d(\Omega)}$ for sufficiently large $n$, so it follows that uniformly $(1 - \hat{\varepsilon}_{a_n}^*/\varepsilon_n)_+$ tends to 0 in probability. This concludes the proof of Theorem 6.1.

## REFERENCES

[1] BENJAMINI, Y. and HOCHBERG, Y. (1995). Controlling the false discovery rate: A practical and powerful approach to multiple testing. *J. Roy. Statist. Soc. Ser. B* **57** 289–300. MR1325392

[2] CAI, T., JIN, J. and LOW, M. G. (2006). Estimation and confidence sets for sparse normal mixtures. Technical report, Dept. Statistics, The Wharton School, Univ. Pennsylvania. Available at www.arxiv.org/abs/math/0612623.

[3] CAI, T. and LOW, M. G. (2004). An adaptation theory for nonparametric confidence intervals. *Ann. Statist.* **32** 1805–1840. MR2102494

[4] Donoho, D. (1988). One-sided inference about functionals of a density. *Ann. Statist.* **16** 1390–1420. MR0964930

[5] Donoho, D. and Jin, J. (2004). Higher criticism for detecting sparse heterogeneous mixtures. *Ann. Statist.* **32** 962–994. MR2065195

[6] Efron, B. (2004). Large-scale simultaneous hypothesis testing: The choice of a null hypothesis. *J. Amer. Statist. Assoc.* **99** 96–104. MR2054289

[7] Genovese, C. and Wasserman, L. (2004). A stochastic process approach to false discovery control. *Ann. Statist.* **32** 1035–1061. MR2065197

[8] Ingster, Y. I. (1999). Minimax detection of a signal for $l_n^p$-balls. *Math. Methods Statist.* **7** 401–428. MR1680087

[9] Jin, J. (2004). Detecting a target in very noisy data from multiple looks. In *A Festschrift for Herman Rubin* (A. DasGupta, ed.) 255–286. IMS, Beachwood, OH. MR2126903

[10] Jin, J. (2006). Proportion of nonzero normal means: Universal oracle equivalences and uniformly consistent estimations. Technical report, Dept. Statistics, Purdue Univ.

[11] Jin, J., Peng, J. and Wang, P. (2007). Estimating the proportion of non-null effects, with applications to CGH lung cancer data. Working manuscript.

[12] Le Cam, L. and Yang, G. L. (1990). *Asymptotics in Statistics*: *Some Basic Concepts.* Springer, New York. MR1066869

[13] Maraganore, D. M., de Andrade, M. et al. (2005). High-resolution whole-genome association study of Parkinson disease. *Amer. J. Human Genetics* **77** 685–693.

[14] Meinshausen, N. and Rice, J. (2006). Estimating the proportion of false null hypotheses among a large number of independently tested hypotheses. *Ann. Statist.* **34** 373–393. MR2275246

[15] Shorack, G. R. and Wellner, J. A. (1986). *Empirical Processes with Applications to Statistics.* Wiley, New York. MR0838963

T. T. Cai
M. G. Low
Department of Statistics
The Wharton School
University of Pennsylvania
Philadelphia, Pennsylvania 19104-6340
USA
E-mail: tcai@wharton.upenn.edu
        lowm@wharton.upenn.edu

J. Jin
Department of Statistics
Purdue University
West Lafayette, Indiana 47907
USA
E-mail: jinj@stat.purdue.edu